\documentclass[11pt]{article}
\usepackage[dvips]{graphicx}
\usepackage{textcomp}
\usepackage{amsbsy}
\usepackage{latexsym}
\usepackage[mathscr]{eucal}
\usepackage{amsfonts}
\usepackage{amssymb}
\usepackage{mathtools}
\usepackage[usenames]{color}
\usepackage[dvipsnames]{xcolor}
\usepackage{comment}

\usepackage{hyperref}

\usepackage{fullpage}
\newtheorem{lemma}{Lemma}[section]
\newtheorem{prop}[lemma]{Proposition}

\newtheorem{theorem}[lemma]{Theorem}
\newtheorem{remark}[lemma]{Remark}

\def\beginproof{\noindent{\bf Proof:}~ }
\def\endproof{\hfill\rule{1.5mm}{1.5mm}\\[2mm]}

\def\cV{{\cal V}}

\newcommand{\be}{\begin{equation}}
\newcommand{\ee}{\end{equation}}

\def\anew{\color{red} }

\def\RR{\rm \hbox{I\kern-.2em\hbox{R}}}
\def\NN{\rm \hbox{I\kern-.2em\hbox{N}}}
\def\ZZ{\rm {{\rm Z}\kern-.28em{\rm Z}}}
\def\CC{\rm \hbox{C\kern -.5em {\raise .32ex \hbox{$\scriptscriptstyle|$}}\kern
-.22em{\raise .6ex \hbox{$\scriptscriptstyle |$}}\kern .4em}}
\def\vp{\varphi}
\def\<{\langle}
\def\>{\rangle}
\def\t{\tilde}

\def\e{\varepsilon}

\def\o{\overline}

\def\cT{{\cal T}}

\def\cV{{\cal V}}

\def\cF{{\cal F}}

\def\Chi{\raise .3ex
\hbox{\large $\chi$}} \def\vp{\varphi}
\def\lsima{\hbox{\kern -.6em\raisebox{-1ex}{$~\stackrel{\textstyle<}{\sim}~$}}\kern -.4em}
\def\lsim{\hbox{\kern -.2em\raisebox{-1ex}{$~\stackrel{\textstyle<}{\sim}~$}}\kern -.2em}
\def\[{\Bigl [}
\def\]{\Bigr ]}
\def\({\Bigl (}
\def\){\Bigr )}
\def\[{\Bigl [}
\def\]{\Bigr ]}
\def\({\Bigl (}
\def\){\Bigr )}

\def\R{\mathbb{R}}
\def\E{\mathbb{E}}

\def\T{{\relax\ifmmode I\!\!\hspace{-1pt}T\else$I\!\!\hspace{-1pt}T$\fi}}

\def\P{\mathbb{P}}

\def\lsim{\raisebox{-1ex}{$~\stackrel{\textstyle<}{\sim}~$}}

  \def\NN{N}                  

\def\t#1{\tilde{#1}}

\def\argmin{\mathop{\rm argmin}}

\def\cS{{\cal F}}

\def\bI{{\bf I}}

\def\cT{{\cal T}}

\def\cV{{\cal V}}
\def\cW{{\cal W}}

\def\cS{{\cal S}}

\def\bG{{\bf G}}

\def\bv{{\bf v}}

\def\argmin{\mathop{\rm argmin}}

\newcommand{\bd}{{\bf d}}

\newcommand{\bea}{$$ \begin{array}{lll}}
\newcommand{\eea}{\end{array} $$}

\def \exp{\mathop{\rm    exp}}

\newcommand{\beqn}{\begin{equation}}
\newcommand{\eeqn}{\end{equation}}

\makeatletter
\@addtoreset{equation}{section}
\makeatother

\newcommand\abs[1]{\left|#1\right|}

\begin{document}
\title
{
Near-optimal approximation methods for elliptic PDEs\\
with lognormal coefficients
}
\author{ 
Albert Cohen and Giovanni Migliorati
\thanks{%
   This research was supported by the   
   Institut Universitaire de France; the ERC Adv grant BREAD.   }  }

\maketitle
\date{}

\begin{abstract}
\noindent 
This paper studies numerical methods for the approximation of elliptic PDEs with lognormal coefficients of the form $-{\rm div}(a\nabla u)=f$ where $a=\exp(b)$ and $b$ is a Gaussian random field. The approximant of the solution $u$ is an $n$-term polynomial expansion in the scalar Gaussian random variables that parametrize $b$. We present a general convergence analysis of weighted least-squares approximants for smooth and arbitrarily rough random field, using a suitable random design, for which we prove optimality in the following sense: their convergence rate matches exactly or closely the rate that has been established in \cite{BCDM} for best $n$-term approximation by Hermite polynomials, under the same minimial assumptions on the Gaussian random field. This is in contrast with the current state of the art results for the stochastic Galerkin method that suffers the lack of coercivity due to the lognormal nature of the diffusion field. Numerical tests with $b$ as the Brownian bridge confirm our theoretical findings. 
\end{abstract}

\section{Introduction}

\subsection{Elliptic PDEs with lognormal coefficients}
\noindent
The theoretical and numerical treatment of parametric PDEs with a large number of deterministic or stochastic parameters has been the object of intensive investigation in recent years, see \cite{CD} and references therein. One principal objective is to design numerical approximation
methods that are immune to the curse of dimensionality in the sense that 
approximation rates can be established even in the case of countably many variables.

The present work is devoted to design and study concrete numerical methods that provably
achieve this goal in the case of elliptic PDEs with lognormal coefficients. 
While our analysis can be applied to other models, we focus on the main following prototype:
given a bounded Lipschitz domain $D\subset \mathbb{R}^d$ with $d=1,2,3$ and a function $f\in L^2(D)$,  
we consider the elliptic partial differential equation    
\begin{equation}
-\textrm{div}(a\nabla u)=f, \textrm{ in } D,
\label{eq:model_pde}
\end{equation}
completed with homogeneous Dirichlet boundary conditions: 
\begin{equation}
u=0, \textrm{ on } \partial D. 
\label{eq:bc}
\end{equation}
We assume that the diffusion coefficient $a$ is a random function defined as
\be
\label{eq:diffusion_coef}
a=\exp(b),
\ee
where $b$ is a centered Gaussian process defined on $D$, with covariance function
\be
\label{eq:covariance_function}
C_b(x,x'):= \mathbb{E}(b(x)b(x')), \quad x,x'\in D.
\ee
The motivation for considering such a diffusion model is for instance the description of 
sedimental media with high contrast in ground water modeling.

We are interested in parametric expansions of the form $a=a(y)=\exp(b(y))$ with
\begin{equation}
 b(y)=\sum_{j\geq 1} y_j \psi_j, 
\label{param}
\end{equation}
where $y=(y_j)_{j\geq 1}$ is a sequence of independent standard 
Gaussian random variables and the functions $\psi_j:D\to\mathbb{R}$ are given. 
An instance of the expansion \eqref{param} is 
 the Karhunen-Lo\`eve (KL) expansion
$b=\sum_{j \geq 1} \xi_j \vp_j$ using the $L^2(D)$-orthonormal eigenfunctions $(\varphi_j)_{j\geq 1}$ 
of the integral operator with kernel $C_b$ and defining $\psi_j:=\sqrt{\lambda_j} \vp_j$ with $\lambda_j:=\E(\xi_j^2)$.
Other expansions may be considered as well, possibly more relevant and/or advantageous than KL for our numerical purposes, as recently shown in \cite{BCM}.

The sequence $y$ ranges over the unbounded domain
\begin{equation*}
U:=\mathbb{R}^\mathbb{N}, 
\end{equation*}
sometimes referred as the parameter set. 
In this paper we work with the usual product measure space given by 
\begin{equation}
(U, \mathcal{B}(U), d\gamma)=(\mathbb{R}^{\mathbb{N}}, \mathcal{B}(\mathbb{R}^{\mathbb{N}}), d\gamma),
\label{measure}
\end{equation}
where $\mathcal{B}(U)=\mathcal{B}(\mathbb{R}^{\mathbb{N}})$ denotes
the $\Sigma$-algebra generated by the Borel cylinders
and $d\gamma$ the tensorized Gaussian probability measure on $U$. 

For any given $y\in U$ such that $b(y)\in L^\infty$, Lax-Milgram theory allows us to define
the solution $u(y)$ in the space $V:=H^1_0(D)$ through the variational formulation
\be
\label{var}
\int_D a(y)\nabla u (y)\nabla v=\int_D fv, \quad v \in V.
\ee
Note that for each $y$ the solution $u(y)$ is a function $x\mapsto u(x,y)$ of the spatial variable $x\in D$,
and the gradient $\nabla$ and integration over $D$ in the above formulation are with respect to this variable, which we often omit
to lighten the notation. The standard Lax-Milgram a-priori estimate tell us that
\be
\|u(y)\|_V\leq C \exp(\|b(y)\|_{L^\infty}), 
\label{apriori}
\ee
where $C:=\|f\|_{V'} =\|f\|_{H^{-1}}\leq C_P \|f\|_{L^2}$ and $C_P$ is the Poincar\'e constant of $D$. 
The variational formulation extends to any $f\in V'$, replacing the right side integral by the duality bracket $\<f,v\>$.

\subsection{Approximation results}
\noindent 
Let us recall some results from \cite{BCDM,CM2016b,HS2014} on the analysis and approximation of the 
solution map
\begin{equation*}
y\mapsto u(y),
\end{equation*} 
induced by the model
\eqref{eq:model_pde}--
\eqref{eq:covariance_function}. Such results involve the size properties of the functions $\psi_j$. 
More precisely, in \cite{BCDM,CM2016b} these size properties are described in the mildest form through conditions of the type
\be
\sup_{x\in D}\sum_{j\geq 1} \rho_j |\psi_j(x)| =K <\infty,
\label{rhocond}
\ee
where $\rho=(\rho_j)_{j\geq 1}$ is a given sequence of positive number with prescribed growth. This type of condition has the merit of taking into account the support properties of the functions $\psi_j$,
in contrast to conditions on the summability of the $\|\psi_j\|_{L^\infty}$ which were initially proposed in 
\cite{HS2014}.

The solution map is only defined on the set 
\begin{equation*}
U_0:=\{y\in U \; : \; b(y)\in L^\infty \}.
\end{equation*}
Note that this set is equal to $U$ in the case of finitely many variable, that is, when $\psi_j=0$ for $j>J$.
However, in both finitely or countably many variable cases, the solution map is not uniformly bounded.
We thus study this map in the Bochner spaces
\begin{equation*}
\cV_k=L^k(U,V,d\gamma)=L^k(U,d\gamma)\otimes V,
\end{equation*}
of functions from $U$ to $V:=H^1_0(D)$ 
which are strongly measurable and such that $y\mapsto \|u(y)\|_V$ belongs to
$L^k(U,d\gamma)$.  

One sufficient condition established in \cite{Cha} for membership in $\cV_k$ for all finite $k$ 
is that the covariance function 
$C_b$ has $C^\e$ H\"older smoothness for some $\e>0$.
The following result from \cite{BCDM} gives an alternative criterion in terms
of a very mild growth conditions on the sequence $(\rho_j)_{j\geq 1}$ such that \eqref{rhocond} holds.

\begin{theorem}
\label{meas}
Assume that there exists a sequence $(\rho_j)_{j\geq 1}$ of positive numbers such that \eqref{rhocond} holds
and $\sum_{j\geq 1} \exp(-\rho_j^2)<\infty$. 
Then $U_0$ has full measure
and $\E(\exp(k\|b(y)\|_{L^\infty}))<\infty$ for all $k<\infty$. In turn
the map $y\mapsto u(y)$ belongs to 
$\cV_k$ for all $k<\infty$.
\end{theorem}

Under the above assumptions, the solution map is in particular in $\cV_2$ and
may therefore be expanded in Hermite series 
\begin{equation}
u=\sum_{\nu\in{\cal F}} u_\nu H_\nu(y), \quad H_\nu(y)=\prod_{j\geq 1}H_{\nu_j}(y_j),\quad u_\nu:=\int_U u(y)H_\nu(y) \, d\gamma,
\label{eq:hermite}
\end{equation}
where ${\cal F}$ denotes the set of finitely-supported sequences of nonnegative integers, 
and $(H_k)_{k\geq 0}$ is the sequence of univariate Hermite polynomials. These polynomials are defined by
\begin{equation*}
H_k(t):= \dfrac{(-1)^k}{\sqrt{k!} g(t)} \dfrac{ d^k   }{ dt^k } \left( g(t)\right), 
\qquad g(t)=\dfrac{1}{\sqrt{2\pi}} \exp(-t^2/2),
\end{equation*}
and normalized according to 
$\int_{\mathbb{R}} | H_k(t)|^2 g(t)dt=1$. Note that $(H_\nu)_{\nu \in {\cal F}}$ is an orthonormal basis of $L^2(U,d\gamma)$
and that the Hermite coefficients $u_\nu$ are elements of $V$. 

Truncation of the Hermite series to a finite set $\Lambda\subset \cF$ yield a approximation
$u_\Lambda=\sum_{\nu\in\Lambda} u_\nu H_\nu$
that belongs to the $V$-valued polynomial space
\begin{equation*}
\cV_\Lambda:=V\otimes \mathbb{P}_{\Lambda},
\end{equation*}
where $\mathbb{P}_\Lambda:=\textrm{span}\left\{ H_\nu : \nu \in \Lambda \right\}$. By Parseval's equality
\begin{equation*}
\|u-u_\Lambda\|_{\cV_2}=\(\sum_{\nu\notin \Lambda} \|u_\nu\|_V^2\)^{1/2}.
\end{equation*}
For any given $n$, the set $\Lambda$ of cardinality
$\#(\Lambda)=n$ that minimizes the above tail bound is the
set $\Lambda_n^*$ of indices that corresponds to the $n$ largest
$\|u_\nu\|_V$. This set is not necessarily unique and it is not
accessible in practice. It is however possible to easily access a family of
suboptimal sets $(\Lambda_n)_{n\geq 1}$ for which 
$\|u-u_{\Lambda_n}\|_{\cV_2}$ satifies a decay estimate of the form $n^{-s}$
under certain growth conditions on the weights $(\rho_j)_{j\geq 1}$
such that \eqref{rhocond}Ê holds. 

This is based on two main results from \cite{BCDM}. The first result establishes
a weighted $\ell^2$ estimate for the $\|u_\nu\|_V$.

\begin{theorem}
\label{summability_theo_lognormal}
Let $r\geq 1$ be an integer.
Assume that there exists a positive sequence $\rho=(\rho_j)_{j\geq 1}$ such that 
$\sum_{j\geq 1} \exp(-\rho_j^2)<\infty$ and such that 
\begin{equation}
\sum_{j\geq 1} \rho_j|\psi_j(x)|=K< C_r:=\frac{\ln 2}{\sqrt r}, \quad x\in D.
\label{eq:psicr}
\end{equation}
Then, one has 
\begin{equation}
\label{eq:wl2h}
\sum_{\nu\in {\cal F}} \xi_\nu \|u_\nu\|_V^2<\infty,
\end{equation}
where
\begin{equation}
\xi_\nu:=\sum_{\|\widetilde \nu\|_{\ell^\infty}\leq r}{\nu\choose \widetilde \nu} \rho^{2\widetilde \nu}=\prod_{j\geq 1}\left(\sum_{l=0}^r {\nu_j\choose l}\rho_j^{2l}\right), \quad \quad \rho^\nu:=\prod_{j\geq 1}\rho_j^{\nu_j}, \quad {\nu\choose \widetilde \nu}:=\prod_{j\geq 1}{\nu_j\choose \widetilde \nu_j},
\label{eq:weights_lognormal}
\end{equation}
with the convention that ${k \choose l}=0$ when $l>k$.
The constant bounding this sum depends on $\|f\|_{V'}$, $\sum_{j\geq 1} \exp(-\rho_j^2)$ and on the difference $C_r-K$.
\end{theorem}

The truncation strategy now consists in taking for
$\Lambda_n$ the set corresponding to the $n$ smallest weights
$\xi_\nu$, so that 
\be
\|u-u_{\Lambda_n}\|_{\cV_2}\leq C\xi_{n+1}^{-1/2},
\label{xin}
\ee
where $C$ is the value of the weighted sum in \eqref{eq:wl2h}
and $(\xi_k)_{k\geq 1}$ is the increasing rearrangement of the $(\xi_\nu)_{\nu\in\cF}$.
In order to further bound the right-hand side in \eqref{xin}, we use a
second result that relates the $\ell^q$ summability of the $\rho_j^{-1}$
and $\xi_\nu^{-1/2}$ where the $\xi_\nu$ are the weights in 
\eqref{eq:weights_lognormal}. 

\begin{theorem}
\label{thm:lemmaqmh}
For any $0<q<\infty$ and any integer $r\geq 1$ such that $q> \frac 2 r$, we have
\begin{equation*}
(\rho_j^{-1})_{j\geq 1} \in \ell^q(\mathbb{N})\implies (
\xi_\nu^{-1/2}
)_{\nu\in{\cal F}} \in \ell^q
({\cal F}). 
\end{equation*}
\end{theorem}

Since for any positive decreasing sequence $\alpha=(\alpha_k)_{k\geq 1}\in \ell^q$, one has 
$\alpha_k\leq Ck^{-1/q}$ with $C=\|\alpha\|_{\ell^q}^{1/q}$, the combination of the two above
results yields the following convergence bound.

\begin{theorem}
\label{thm:theoherm}
If \eqref{eq:psicr} holds with $(\rho_j^{-1})_{j\geq 1}\in \ell^q(\mathbb{N})$ for some $0<q<\infty$, then
\begin{equation*}
\|u-u_{\Lambda_n}\|_{\cV_2}\leq Cn^{-s}, \quad s:=\frac 1 q,
\end{equation*}
where $\Lambda_n$ is a set corresponding to the $n$ largest $\xi_\nu^{-1/2}$, with $\xi_\nu$ given by~\eqref{eq:weights_lognormal}.
\end{theorem}

\begin{remark}
Since renormalization of the sequence $\rho$ 
does not affect the fact that $(\rho_j^{-1})_{j\geq 1}\in \ell^q(\mathbb{N})$, 
we may replace \eqref{eq:psicr} by the condition \eqref{rhocond}
without a specific constraint on $K$. However the construction of the 
set $\Lambda_n$ requires to define the $\xi_\nu$ using a normalization of $\rho_j$
such that \eqref{eq:psicr} holds.
\end{remark}

\begin{remark}
\label{remspace}
One of the advantage of the above results
is that condition \eqref{rhocond} takes into account the spatial 
localization properties of the functions $\psi_j$, in contrast to 
the analysis firstly proposed in \cite{HS2014} 
that uses conditions on the norms $\|\psi_j \|_{L^\infty(D)}$. 
This is particularly effective in the case of Gaussian fields $b$
that can be described by multiresolution
expansions of wavelet type, such as the Brownian bridge and the general
class of multivariate Mat\'ern processes as discussed in \cite{BCM}.
\end{remark}

One elementary yet important observation, see Lemma 4 from \cite{CM2016b}, is
that for any integer $r\geq 1$, the sequence $(\xi_\nu)_{\nu \in {\cal F}}$ is monotone increasing,
that is 
\begin{equation*}
\t \nu \leq \nu \implies \xi_{\t \nu} \leq \xi_{\nu},
\end{equation*}
where $\t \nu \leq \nu$ means that $\t \nu_j \leq \nu_j$ for all $j$. This implies that the
$\Lambda_n$ may be chosen so to have the following property: a set $\Lambda$ is
{\em downward closed} if
\begin{equation*}
\nu \in \Lambda \textrm{ and } \t \nu \leq \nu
\implies \t \nu \in \Lambda.
\end{equation*}
This property plays an important role in the 
development of numerical methods for parametric and stochastic PDEs 
like \eqref{eq:model_pde}, and we use it for those developed in
the present paper. We say that the polynomial space 
$\mathbb{P}_\Lambda:=\textrm{span}\left\{ H_\nu : \nu \in \Lambda \right\}$ is downward closed 
whenever the set $\Lambda$ is downward closed. 
Notice that, in this setting, the $H_\nu$ can be replaced by the monomials $y\mapsto y^\nu$ or more generally any tensorized polynomial basis $G_\nu$.
 
\subsection{Numerical methods}
\noindent
Once given the sets $\Lambda_n$, the calculation of the truncated Hermite series $u_{\Lambda_n}$ requires the coefficients $u_\nu$ for all $\nu \in \Lambda_n$. 
These coefficients belong to the space $V$, and in general they cannot be computed exactly. 
The main purpose of the present paper is the development and analysis of numerical methods for the approximation of these coefficients. We aim at methods that provably achieve the convergence rate that can be established for the ideal approximation $u_{\Lambda_n}$, under the same assumption \eqref{rhocond} with $(\rho_j^{-1})_{j\geq 1} \in \ell^q$, and at reasonable computational cost.

One prominent numerical method for parametric elliptic problem is the stochastic Galerkin method,
which is based on a variational formulation 
in both space and parametric variable. In the present case, the most natural
such formulation is obtained by integrating the spatial variational formulation \eqref{var} on the 
parametric domain $U$ with respect to the measure $\gamma$: find $u$ in $\cV_2$ such that 
\begin{align*}
A(u,v)=L(v), \quad v\in \cV_2,
\end{align*}
where the bilinear and linear forms $A$ and $L$ are defined by
\begin{align*}
&A(u,v):=\int_U \(\int_D a(y) \nabla u(y)\nabla v(y)\) d\gamma(y)=\int_{U\times D} a(x,y)\nabla u(x,y).\nabla v(x,y)\, dx \,d\gamma(y),\\
&L(v):=\int_{U}\(\int_D fv(y)\)d\gamma(y)=\int_{U\times D}f(x)v(x,y) \,dx\, d\gamma(y). 
\end{align*}

The Galerkin method then amounts to restricting the above variational formulation to a subspace
of $\cV_2$ such as the polynomial space $\cV_{\Lambda_n}$. 
However this approach leads to a generally ill-posed problem, due to the 
degeneracy of the lognormal diffusion coefficient $a(x,y)$ which is positive but may 
attain arbitrarily small or large values. As a consequence the bilinear form $A$
is neither continuous nor coercive on $\cV_2$, 
thus preventing the use of Lax-Milgram
theory and Cea's lemma. More precisely, the quadratic form associated to the bilinear form $A$ is given by
\begin{equation*}
A(u,u)=\int_U \|u(y)\|_{a(y)}^2 d\gamma(y),
\end{equation*}
where 
\begin{equation*}
\|v\|_{a(y)}^2 := \int_D a(y) |\nabla v|^2.
\end{equation*}
This $y$-dependent norm and the fixed $V$-norm satisfy the equivalence
\begin{align}\label{equivaV}
\exp\left(-\|b(y)\|_{L^\infty}\right)\|v\|^2_V\leq \|v\|^2_{a(y)} \leq \exp\left(+\|b(y)\|_{L^\infty}\right)\|v\|^2_V, \quad v\in V,
\end{align}
with sharp constants that depend on $y$ and are not uniformly controlled over $U$. Therefore $A(u,u)$
is not equivalent to $\|u\|_{\cV_2}^2$.

Several alternatives have been developed in the literature in order 
to circumvent this obstacle. One approach first proposed in \cite{GS2009}
is to modify the choice of spaces in the formulation and establish proper 
inf-sup conditions. Another approach introduced in \cite{Gitt10}
is based on the use of an auxilliary Gaussian measure $d\tilde{\gamma}$ 
and corresponding formulation
\begin{align*}
\tilde{A}(u,v)=\tilde{L}(v),
\end{align*}
where $\tilde{A}$ and $\tilde{L}$ are defined similarly as $A$ and $L$ with $\gamma$ replaced by $\t \gamma$.
This second approach is used in \cite{HS2014} where approximation rates for the Galerkin method on appropriate polynomial spaces are established. This however requires conditions of the type $(j\|\psi\|_{L^\infty})_{j\geq 1}\in \ell^p$ for some $p\leq 1$.

It is not clear if the methods in \cite{GS2009,Gitt10,HS2014} may yield the same convergence rate that can be established for the ideal approximation $u_{\Lambda_n}$, under the mild assumption \eqref{rhocond} with $(\rho_j^{-1})_{j\geq 1} \in \ell^q$. 
In the present paper we adopt a different strategy based on weighted least squares that allows us to achieve this goal by computationally feasible methods. 

\subsection{Outline of the paper}
\noindent
As a first step, discussed in \S 2, we proceed to a discretization in the spatial variable $x$, by replacing for each $y\in U$ the solution $u(y)\in V$ defined in \eqref{var} by its discrete counterpart $u_h(y)\in V_h$ defined by the Galerkin method
\begin{equation*}
\int_D a(y) \nabla u_h(y)\nabla v_h=\int_D f v_h, \quad v_h\in V_h.
\end{equation*}
Here $V_h\subset V$ is a finite-dimensional approximation space, such as a finite element space.
Our main result in this section gives convergence bounds as $h\to 0$ for
the $\cV_2$ error between the solution map $y\mapsto u(y)$ and its discrete 
counterpart $y\mapsto u_h(y)$. These bounds are established by classical elliptic
regularity estimates when the random fields $b$ and $a$ are sufficiently smooth,
and extended to arbitrarily rough fields by interpolation theory arguments.
We illustrate these bounds for Gaussian fields $b$ described by multiresolution
expansions, already mentionned in Remark \ref{remspace}.

In \S\ref{sec:wls}, we discuss the computation of the polynomial
approximation by the discrete weighted least-squares
method: we compute particular snapshots $u(y^i)$ for $i=1,\dots,m$, up
to space discretization in $V_h$, where $y^i\in U$ are randomly sampled according to some measure $d\mu$,
and obtain the approximation 
by minimizing a criterion of the form
\begin{equation*}
\sum_{i=1}^m w(y^i) \|u(y^i)-v(y^i)\|_V^2,
\end{equation*}
over $v \in \cV_{\Lambda_n}$. Here, we exploit the results from \cite{CM2016}, which show that a proper choice of the weight function $w$ and sampling measure $d\mu$ provides optimal error estimates in expectation, with $m$ of the same order as $n$ up to a logarithmic factor. 

In \S\ref{ref:examples_random_field}, we discuss the convergence rates guaranteed by our error analysis, for particular 
Gaussian processes $b$ of practical interest: univariate Brownian processes
and multivariate Mat\'ern covariances. The chosen expansion \eqref{param} is of wavelet
type, building up on the results in \cite{BCM}. With such choices, it is shown that the
obtained rates for spatial and parametric approximation are the same and 
proportional to the smoothness parameter $\delta>0$ appearing in the covariance
function for Mat\'ern processes. This is in contrast with the use of standard Karhunen-Loeve basis
for which the provable parametric rate is lesser than the spatial rate.

Finally numerical tests are presented in \S\ref{sec:num_illust} for the weighted least-squares method, for the Brownian processes truncated up to various parametric dimensions.

\section{Space discretization}
\label{sec:SD}
\noindent
In this section, we introduce and study a semi-discretization of the solution $u$ to the parametric elliptic problem \eqref{eq:model_pde}, with respect to the spatial variable $x\in D$. For this purpose, we denote by $(V_h)_{h>0}$ a family of 
closed finite-dimensional subspaces of $V$ of dimension 
\begin{equation*}
n_h:=\dim(V_h),\quad  {\rm where}\quad n_h\to +\infty \quad {\rm as}\quad h\to 0.
\end{equation*}
For any $y\in U$ such that $b(y)\in L^\infty(D)$, we define $u_h(y)\in V_h$ as the Galerkin approximation to $u(y)$:
\begin{align}
\label{SpaceDiscrPara}
\int_D a(x,y)\nabla u_h (x,y)\nabla v_hdx=\int_D fv_hdx, \quad v_h\in V_h.
\end{align}
Note that for a given $y$, each query of the approximate solution map
\begin{equation*}
y\mapsto u_h(y),
\end{equation*}
amounts to solving an $n_h\times n_h$ linear system.
Similar to the exact solution map $y\mapsto u(y)$, the approximate solution map
satisfies the following measurability and integrability property under mild growth
conditions on the sequence $(\rho_j)_{j\geq 1}$ such that \eqref{rhocond} holds.

\begin{prop}
\label{hmeas}
If the assumptions of Theorem \ref{meas}Ê hold, then
the map $y\mapsto u_h(y)$ is measurable and belongs to 
$\cV_k$ for all $k<\infty$.
\end{prop}
\noindent
By the previous proposition, the approximate solution map $y\mapsto u_h(y)$ is in $\cV_2$ so that we can expand it into multidimensional Hermite series
\begin{equation*}
u_h(y)=\sum_{\nu\in\cF} u_{\nu,h} H_\nu(y),
\end{equation*}
where the coefficients $u_{\nu,h}$ belong to $V_h$.  Let $\Lambda$ be a given subset of $\mathcal{F}$ with finite cardinality $n$. 
We introduce the subspace 
\begin{align*}
\cV_{\Lambda,h}:=V_h\otimes \mathbb{P}_\Lambda 
\end{align*}
of $\cV_2$, which has dimension $n_h\times n$. 
In \S\ref{sec:wls} 
we construct
an approximation $u_{\Lambda,h}\in V_{\Lambda,h}$ to the solution map 
by the weighted least-square method, 
which can be viewed as a vehicle to obtain discretizations in the parametric variable.
The error estimate in $\cV_2$ between $u$ and $u_{\Lambda,h}$ may be bounded according to 
\begin{align}\label{FullError}
\|u-u_{\Lambda,h}\|_{\cV_2}\leq \|u-u_{h}\|_{\cV_2}+\|u_h-u_{\Lambda,h}\|_{\cV_2},
\end{align}
where the first and second terms on the right-hand side represents the spatial and parametric discretization errors, respectively.
The estimation of the second term depends on the chosen discretization method
in the parametric variable, and it is discussed in \S\ref{sec:wls}. Here, we concentrate on the first
error term $\|u-u_{h}\|_{\cV_2}$.  The rate of approximation 
which is achievable for this error term depends on the  smoothness of the Gaussian field $b$.

\subsection{Smooth Gaussian field}\label{SGF}
Let us introduce the spatial regularity space
\begin{equation*}
W:=\{v\in V\; : \; \Delta v\in L^2(D)\},
\end{equation*}
equipped with the norm 
\begin{equation*}
\|v\|_W:=\|\Delta v\|_{L^2}.
\end{equation*}
We assume that the spaces $V_h$ possess
the following approximation property with respect to 
the space $W$: there exist $C>0$ and $t>0$ such that
\be
\label{ApproxProp}
 \underset{v_h\in V_h}{\min}\|v-v_h\|_V\leq C n_h^{-t} \|v\|_W,\quad  v\in W, \quad h>0.
\ee
It is well known that $W$ coincides with $V\cap H^2(D)$ with equivalent norms
when the domain $D$ is either convex or has $C^{1,1}$ boundary.
Therefore, in this case, taking for $(V_h)_{h>0}$ the $\P_1$ Lagrange finite elements
associated to a regular family of quasi-uniform conforming simplicial partitions $(\cT_h)_{h>0}$ where $h$ is
the mesh-size of $\cT_h$, we may use results from classical
finite-element approximation theory \cite{Cia} 
\be
\label{ApproxProph}
 \underset{v_h\in V_h}{\min}\|v-v_h\|_V\leq C h |v|_{H^2},\quad  v\in H^2(D), \quad h>0,
\ee
to conclude that \eqref{ApproxProp} hold with rate $t=\frac 1 d$. It is also well known that this
rate can be maintained in the case where $D$ is a nonconvex polygon or polyhedron provided
that the spaces $V_h$ are enriched by appropriate local refinement near the re-entrant corner
or edges.

In order to study the error term $\|u-u_{h}\|_{\cV_2}$, we need to study the integrability of
$y\mapsto u(y)$ in the sense of the Bochner space 
\begin{equation*}
\cW_k=L^k(U,W,\gamma),
\end{equation*}
which is obviously smaller than $\cV_k$. We begin by observing that for the nonparametric equation \eqref{eq:model_pde}, the solution has $W$-regularity when in addition to being elliptic, the coefficient $a$ belongs to $W^{1,\infty}(D)$. Indeed, the equation then writes
\begin{equation*}
\Delta u = a^{-1}(\nabla a\cdot \nabla u-f),
\end{equation*}
and since it was assumed that $f\in L^2(D)$, one has
\be\label{ineq:AprioriW}
\|u\|_W \leq a_{\min}^{-1} (\|\nabla a\|_{L^\infty} \|u\|_V+\|f\|_{L^2})\leq 
a_{\min}^{-1} (\|\nabla a\|_{L^\infty} a_{\min}^{-1}\|f\|_{V'}+\|f\|_{L^2}),
\ee
where $a_{\min}=\min_{x\in D} a(x)>0$.
If $a=\exp(b)$ for some $b\in L^\infty$, this gives
\be
\|u\|_W \leq  \exp(\|b\|_{L^\infty})(\|\nabla b\|_{L^\infty}\exp(2\|b\|_{L^\infty}) \|f\|_{V'}+\|f\|_{L^2}).
\label{West}
\ee
Note that $\|f\|_{V'}\leq C_P\|f\|_{L^2}$ where $C_P$ is the Poincar\'e constant.
This leads to the following result that gives a condition
ensuring $L^k$-integrability of $\|u(y)\|_{W}$ 
with respect to the Gaussian measure $d\gamma$.

\begin{prop}\label{SpaceReg} 
If the Gaussian field $b$ belongs to $W^{1,\infty}(D)$ almost surely
with $\E(\|\nabla b\|_{L^\infty}^k)<\infty$ for all $k<\infty$, then
 $u\in \cW_k$ for all $k<\infty$.
\end{prop}

\beginproof
Using \eqref{West}, we may write for any $y$ such that $b(y)\in W^{1,\infty}(D)$, 
\begin{equation*}
\|u(y)\|_{W}\leq \exp(\| b(y) \|_{L^\infty})\bigg(\|\nabla b(y) \|_{L^\infty}\exp(2\|b(y)\|_{L^\infty} )\|f\|_{V'}+\|f\|_{L^2}\bigg).
\end{equation*}
By application of Cauchy-Schwarz inequality and using Theorem \ref{meas}, we find that $u\in \cW_k$
with
\begin{equation*}
\|u\|_{\cW_k} \leq\E(\|\nabla b\|_{L^\infty}^{2k})^{\frac 1 {2k}}\E(\exp(6k\|b(y)\|_{L^\infty}))^{\frac 1 {2k}} \|f\|_{V'}
+\|f\|_{L^2}\E(\exp(k\|b(y)\|_{L^\infty}))^{\frac 1 k},
\end{equation*}
where the terms on the right are controlled from the assumptions.
Indeed, if the Gaussian field $b$ belongs to $W^{1,\infty}(D)$ almost surely, this implies that $b$ induces a Gaussian measure on the separable Banach space $C^0(\overline{D})$. Then, by Fernique's Theorem (see Theorem 2.8.5 in \cite{B98}), we find that
$\mathbb{E}(\exp (k \|b\|_\infty))<\infty$ for all $k<\infty$. Finally, by assumption, $\E(\|\nabla b\|_{L^\infty}^k)<\infty$ for all $k<\infty$.\endproof

As an immediate consequence, we obtain a convergence rate 
of the spatial discretization error for regular Gaussian fields.

\begin{theorem}\label{SDError}
If the assumptions of Proposition \ref{SpaceReg} hold, then
there exists $C>0$ such that, for any $h>0$, we have
\begin{align*}
\|u-u_{h}\|_{\cV_2}\leq Cn_h^{-t},
\end{align*}
where $t$ is the approximation rate for functions in $W$, as in \eqref{ApproxProp}.
\end{theorem}

\beginproof
Using the norm equivalence (\ref{equivaV}) and the fact
that the Galerkin approximation $u_h(y)$ is the orthogonal projection of 
$u(y)$ onto $V_h$ in the sense of the $a(y)$ norm, we obtain
\begin{align*}
\|u(y)-u_{h}(y)\|^2_V &\leq \exp\left( \|b(y)\|_{L^\infty}\right)\|u(y)-u_{h}(y)\|^2_{a(y)} \\
 &\leq \exp\left(\|b(y)\|_{L^\infty}\right)\underset{v_h\in V_h}{\min}\|u(y)-v_{h}\|^2_{a(y)}\\
 & \leq \exp\left(2\|b(y)\|_{L^\infty}\right)\underset{v_h\in V_h}{\min}\|u(y)-v_{h}\|^2_{V}.
\end{align*}
By the approximation property (\ref{ApproxProp}), this yields
\begin{align*}
\|u(y)-u_{h}(y)\|^2_V\leq C\exp\left(2\|b(y)\|_{L^\infty}\right) n_h^{-2t}\|u(y)\|^2_{W}.
\end{align*}
Therefore, by Cauchy-Schwarz inequality, we find
\begin{align*}
\|u-u_{h}\|^2_{\cV_2}\leq Cn_h^{-2t} \E(\exp(4\|b(y)\|_{L^\infty})^{1/2} \|u\|_{\cW_4}^{2},
\end{align*}
which is the claimed estimate up to a change in the constant $C$.
\endproof

We end this subsection by discussing the validity of the additional
assumption in Proposition \ref{SpaceReg}, namely that the
Gaussian field $b$ belongs to $W^{1,\infty}(D)$ almost surely and $\E(\|\nabla b\|_{L^\infty}^k)<\infty$
for all $k<\infty$. Standard conditions for this to hold can be stated in terms of the smoothness of 
the covariance function $C_b$. Roughly speaking this function should be a bit smoother
than $C^2$ along the diagonal $x=x'$. In the simple case of one-dimensional second order stationary stochastic processes where $C_b(x,x')=c_b(x-x')$, such a result (see Section $9$ in \cite{CrLea67}) reads, for some $a>3$ and $\ell <+\infty$,
\begin{equation*}
c_b(r)=c_b(0)-\frac{\ell r^2}{2}+O\left(\frac{r^2}{|\log(|r|)|^a}\right),
\end{equation*}
when $r$ tends to $0$.

Based on this type of criterion we find that $b$ is almost surely in $C^1(\o D)$ and induces
a Gaussian measure on this separable Banach space. Then Fernique's Theorem (Theorem 2.8.5 in \cite{B98}) 
implies that $\|b(y)\|_{C^1}$ belongs to $L^k(U,d\gamma)$ for all $k<\infty$
and therefore $\E(\|\nabla b\|_{L^\infty}^k)<\infty$ for all $k<\infty$.

Many relevant random fields are however excluded from this analysis since 
their sample path are not regular enough.  This is typically the case for Brownian processes 
and their multivariate extensions
and for the Mat\'ern Gaussian fields with small smoothness parameter,
that are discussed in \S 6.  In the next subsection, 
we refine our analysis in order to include these relevant 
 rough Gaussian fields.

\subsection{Rough Gaussian field}
The previous analysis assumes the fact that the Gaussian field $b$ is smooth in the sense that $\nabla b$ exists almost surely and is in $L^\infty(D)$. However, the Gaussian field $b$ is usually not a regular field and its sample paths possess some H\"{o}lderian regularity with index strictly less than $1$. Thus, we need to perform a discretization in the spatial variable taking into account the roughness of the sample paths of $b$. First of all, we introduce some notation. For $0<\alpha<1$, we denote by $C^{\alpha}(D)$ the Banach space of H\"{o}lder continuous functions whose norm is given by:
\begin{align*}
\|b\|_{C^\alpha}:=\|b\|_{L^\infty}+\underset{x_1,x_2\in D
\atop x_1 \ne x_2}{\sup}\dfrac{|b(x_1)-b(x_2)|}{|x_1-x_2|^{\alpha}}.
\end{align*}
We denote as well by $|\cdot|_{C^\alpha}$ the H\"{o}lderian seminorm defined by:
\begin{align*}
|b|_{C^\alpha}:=\underset{x_1,x_2\in D
\atop x_1\ne x_2}{\sup}\dfrac{|b(x_1)-b(x_2)|}{|x_1-x_2|^{\alpha}}.
\end{align*} 

In the following, we use the real interpolation method in order to find the suitable scale of functional spaces in which the solution to the PDE belongs, in the case of such rough fields, as well as the spatial approximation results for such spaces. We recall briefly the definition of interpolation spaces, see \cite{BL} for a general treatment. For general Banach spaces $X$ and $Y$ which are compatible, we define the $K$-functional of any elements $g$ of $X+Y$ by:
\begin{equation*}
K(g,t,X,Y):=\inf\{\|g_X\|_X+t\|g_Y\|_Y\; : \;  g_X\in X,\, g_Y\in Y,\, g=g_X+g_Y\},\,\quad t\geq 0,
\end{equation*}
For any $0<\theta<1$ and for any $0<q\leq \infty$, we define the interpolation space $[X,Y]_{\theta,q}$ as the set of functions of $X+Y$ such that the following quantity is finite
\begin{align*}
\|g\|_{[X,Y]_{\theta,q}}=\|g\|_{\theta,q}:=\left\{
    \begin{array}{ll}
        \bigg(\int_{0}^{+\infty}(t^{-\theta}K(g,t,X,Y))^q\frac{dt}{t}\bigg)^{\frac{1}{q}}, &  q<\infty,\\
       \underset{t>0}{\sup}\, t^{-\theta}K(g,t,X,Y),  & q=\infty.
    \end{array}
\right.
\end{align*}
In the case where $Y\subset X$ which is of interest for our present discussion, the $K$-functional also writes
\begin{equation*}
K(g,t,X,Y):=\inf\{\|g-g_Y\|_X+t\|g_Y\|_Y \; : \;  g_Y\in Y\},\,\quad t\geq 0.
\end{equation*}
In particular, one has $K(g,t,X,Y)\leq \|g\|_X$ and the above interpolation norm may
be equivalently defined by taking the integral or supremum on $t\in [0,1]$.
We shall use the particular well-known case of interpolation spaces: if $D$ is a bounded Lipschitz domain then
one has
\begin{equation*}
C^\alpha(D)=[L^\infty,W^{1,\infty}]_{\alpha,\infty},
\end{equation*}
with equivalent norms. As a substitute to the estimates 
\eqref{ineq:AprioriW} and \eqref{West} in $W$ for smooth fields, the following result gives estimates 
in the interpolation space
\begin{align*}
W^\alpha:=[V,W]_{\alpha,\infty}.
\end{align*}
for the solution to the elliptic problem \eqref{eq:model_pde} when the fields $a$ and $b$ are in $C^{\alpha}(D)$.

\begin{theorem}\label{theo:Besov}
Assume that $D$ is a bounded Lipschitz domain of $\R^d$, 
that $a\in C^{\alpha}(D)$ with $a_{\min}:=\min_{x\in D} a(x)>0$ and that $f\in L^2(D)$. Then, the solution $u$ to the elliptic problem \eqref{eq:model_pde} belongs to $W^\alpha$ and,
\begin{align}
\label{uwalpha}
\|u\|_{W^\alpha}\leq C \( a_{\min}^{-1}+ \|a\|^{\frac{1}{\alpha}}_{C^\alpha} a_{\min}^{-\frac {1+\alpha}{\alpha} }\ \)\|f\|_{L^2}. 
\end{align}
For $a$ of the form $a=\exp(b)$ with $b\in C^\alpha$, one has
\be
\label{uwalphab}
\|u\|_{W^\alpha}\leq C \left(\exp(\|b\|_{L^{\infty}})+(1+2\|b\|_{C^\alpha})^{\frac{1}{\alpha}}\exp
\left( \frac{2+\alpha}{\alpha}\|b\|_{L^\infty} \right) \right)\|f\|_{L^2}.
\ee
The constants $C$ in the above estimate depend on $\alpha$ and $D$ only.
\end{theorem}

\beginproof
We make use of a standard stability result for second order elliptic PDEs. If $u$ and $\t u$ are two weak solutions
to \eqref{eq:model_pde}
with diffusion coefficients $a$ and $\t a$, respectively, and with the same data $f$, then one has
\be\label{ineq:StabEst}
\|u-\t u\|_V \leq    \|a-\t a\|_{L^\infty}  (a_{\min}\t a_{\min})^{-1}\|f\|_{V'},
\ee
where $a_{\min}$ and $\t a_{\min}$ are the minimal values of $a$ and $\t a$ on $D$. This is easily checked 
by equating the variational formulations for both solutions and taking $v=u-\t u$ as a test function. 

Since $C^\alpha(D)=[L^\infty,W^{1,\infty}]_{\alpha,\infty}$ with equivalent norms, 
for any $\varepsilon>0$ there exists $a_\varepsilon\in W^{1,\infty}(D)$ such that 
\be
\label{ineq:InterHol1}
\|a-a_\varepsilon\|_{L^\infty} \leq C \|a\|_{C^\alpha} \varepsilon^\alpha,
\ee
and 
\be
\label{ineq:InterHol2}
 |a_\varepsilon|_{W^{1,\infty}}\leq C \|a\|_{C^\alpha} \varepsilon^{\alpha-1},
\ee
where the constant $C$ only depends on $\alpha$ and $D$. Moreover, 
if $\e\leq \varepsilon_0=\left(a_{\min}/(2\|a\|_{C^\alpha}  C)\right)^{1/\alpha}$, 
we have
\be\label{ineq:ellreg}
a_\varepsilon(x)\geq \frac{a_{\min}}{2}, \quad x\in D.
\ee
For such values of $\e$, we denote by $u_\varepsilon$ the solution to \eqref{eq:model_pde} with $a_\varepsilon$ as diffusion coefficient. We note that $u_\varepsilon$ belongs to $W$ and satisfies 
\begin{equation*}
\|u-u_\varepsilon\|_V \leq   2a_{\min}^{-2}C \|a\|_{C^\alpha}  \|f\|_{V'}\varepsilon^\alpha,
\end{equation*}
by \eqref{ineq:InterHol2} combined with \eqref{ineq:StabEst}, as well as 
\begin{equation*}
\|u_\varepsilon\|_W \leq  2a_{\min}^{-1} (C \|a\|_{C^\alpha}2a_{\min}^{-1}  \varepsilon^{\alpha-1}\|f\|_{V'}+\|f\|_{L^2}),
\end{equation*}
by \eqref{ineq:InterHol2} combined with \eqref{ineq:AprioriW}. 

We replace $\|f\|_{V'}$ by $\|f\|_{L^2}$ up to a multiplication of $C$ by the Poincar\'e constant,
and use $u_\e$ to estimate $K(u,t,V,W)$ for the relevant range $0<t\leq 1$, by writing
\begin{equation*}
K(u,t,V,W)\leq \( 2C\|a\|_{C^\alpha}a_{\min}^{-2}\varepsilon^\alpha+2t a_{\min}^{-1}(C \|a\|_{C^\alpha}2a_{\min}^{-1} \varepsilon^{\alpha-1}+1) \)\|f\|_{L^2}.
\end{equation*}
Thus, for any $0<\varepsilon<\varepsilon_0\wedge \varepsilon_1$, with $\varepsilon_1=\left(2C\|a\|_{C^\alpha}
a_{\min}^{-1}\right)^{1/(1-\alpha)}$
\begin{equation*}
K(u,t,V,W)\leq  C\( 2\|a\|_{C^\alpha}a_{\min}^{-2}\varepsilon^\alpha+8t a_{\min}^{-2} \|a\|_{C^\alpha}\varepsilon^{\alpha-1} \)\|f\|_{L^2}.
\end{equation*}
For $t\leq 1$, we take $\varepsilon=t(\varepsilon_0\wedge \varepsilon_1)$ in the previous inequality to obtain,
\begin{equation*}
K(u,t,V,W)\leq C\left(2\|a\|_{C^\alpha}a_{\min}^{-2} (\varepsilon_0\wedge \varepsilon_1)^{\alpha}+8 a_{\min}^{-2} \|a\|_{C^\alpha}(\varepsilon_0\wedge \varepsilon_1)^{\alpha-1}\right)t^\alpha\|f\|_{L^2}.
\end{equation*}
Using $\varepsilon_0\wedge \varepsilon_1\leq \varepsilon_0$ and $\varepsilon_0\wedge \varepsilon_1\geq \varepsilon_0\varepsilon_1/(\varepsilon_0+\varepsilon_1)$, we arrive at
\begin{equation*}
t^{-\alpha}K(u,t,V,W)\leq C \( a_{\min}^{-1}+ \|a\|^{\frac{1}{\alpha}}_{C^\alpha} a_{\min}^{-\frac {1+\alpha}{\alpha} }\ \)\|f\|_{L^2},
\end{equation*}
which gives \eqref{uwalpha}.
Finally, when $a=\exp(b)$ with $b\in C^\alpha$, we use $\|a\|_{C^{\alpha}}\leq \exp(\|b\|_{L^\infty})(1+2\|b\|_{C^\alpha})$ together with $a_{\min}^{-1}\leq \exp(\|b\|_{L^\infty})$ to obtain \eqref{uwalphab}.
\endproof

\begin{remark}
The interpolation argument used in the proof of the above result also allows us to treat the case of data $f$ that are rougher than $L^2$, by using the stability estimate
\begin{equation*}
\|u-\t u\|_V \leq    \|a-\t a\|_{L^\infty}  (a_{\min}\t a_{\min})^{-1}\|f\|_{V'}+\frac 1 {a_{\min} }\|f-\t f\|_{V'}.
\end{equation*}
in place of \eqref{ineq:StabEst} when $u$ and $\t u$ are associated with different data $f$ and $\t f$. Using the fact that $H^{-1+\alpha}=[H^{-1},L^2]_{\alpha,2} \subset [H^{-1},L^2]_{\alpha,\infty}$, the same argument shows that the solution $u$ belongs to $W^\alpha$ when $a\in C^\alpha$ and $f\in H^{-1+\alpha}$. In the case where $W$ coincides with $V\cap H^2(D)$ (for example, when $D$ is either convex or has $C^{1,1}$ boundary), we find that
\begin{equation*}
W^\alpha=[V,V\cap H^2(D)]_{\alpha,\infty}=V\cap [H^1(D),H^2(D)]_{\alpha,\infty}= V \cap B^{1+\alpha,2}_{\infty}(D),
\end{equation*}
where $B^{1+\alpha,2}_{\infty}$ is the usual Besov space, which contains 
all Sobolev spaces $H^{1+\beta}$ for $\beta<\alpha$.
The fact that $a \in C^\alpha$ and $f\in H^{-1+\alpha}$ implies that
$u\in H^{1+\beta}$ for $\beta<\alpha$ is proved in \cite{Hac} by different techniques,
see Theorem $9.19$ therein.
\end{remark}
\noindent
In order to obtain an estimate on the quadratic error $\|u-u_{h}\|_{\cV_2}$ similar to the one obtained in Theorem \ref{SDError}, we introduce the Bochner type space
\begin{equation*}
\cW_k^\alpha=L^k(U,W^\alpha,\gamma),
\end{equation*}
for $0<\alpha<1$ and for all $k<\infty$. We now state a result that comes as a substitute to Proposition \ref{SpaceReg} in the case of rough fields. 

\begin{prop}
\label{spaceregrough}
Assume that $b$ is a Gaussian field which belongs to $C^\alpha(D)$ almost surely and such that 
$\mathbb{E}( \|b\|^k_{C^\alpha})<\infty$ for all $k<\infty$. Then, $u\in \cW_k^\alpha$ for all $k<\infty$.
\end{prop}

\beginproof
Based on Theorem \ref{theo:Besov}, for all $y$ such that $b(y)\in C^\alpha(D)$,
\begin{align*}
\|u(y)\|_{W^\alpha}\leq C \(\exp(\|b(y)\|_{L^{\infty}})+(1+2\|b(y)\|_{C^\alpha})^{\frac{1}{\alpha}}\exp\(\frac{2+\alpha}{\alpha}\|b(y)\|_{L^\infty}\) \)\|f\|_{L^2}.
\end{align*}
Since the Gaussian field $b$ induces a Gaussian measure on the separable Banach space $C^0(\overline{D})$ Fernique's theorem (see Theorem $2.8.5$ in \cite{B98}) implies that $\mathbb{E}(\exp(k\|b(y)\|_{L^\infty}))<\infty$ for all $k<\infty$. Moreover, by assumption $\mathbb{E} (\|b\|^k_{C^\alpha})<\infty$ for all $k<\infty$. Then, Cauchy-Schwarz inequality concludes the proof of the proposition.
\endproof
\noindent
Now, we notice that functions in the space $W^{\alpha}$ can be approximated at a certain rate. Indeed, the approximation property \eqref{ApproxProp} for the space $W$
and the sequence $(V_h)_{h>0}$ induces a similar property for the space $W^\alpha$
with the same sequence $(V_h)_{h>0}$, namely
\be
\label{ApproxPropalpha}
 \underset{v_h\in V_h}{\min}\|v-v_h\|_V\leq C n_h^{-\alpha t} \|v\|_{W^\alpha},\quad  v\in W^\alpha, \quad h>0.
\ee
This is readily checked by applying the interpolation inequality 
\begin{equation*}
\|T\|_{[X,Y]_{\alpha,\infty}\to X} \leq \|T\|_{Y\to X}^\alpha \|T\|_{X\to X}^{1-\alpha},
\end{equation*}
with $X=V$ and $Y=W$ to $T=I-P_{V_h}$ where $P_{V_h}$ is the $V$-orthogonal projection onto $V_h$.
This leads to the following approximation result.

\begin{theorem}\label{thm:RoughError}
If the assumptions of Proposition \ref{spaceregrough} hold, then
there exists $C>0$ such that, for any $h>0$, we have
\begin{align*}
\|u-u_{h}\|_{\cV_2}\leq Cn_h^{-\alpha t},
\end{align*}
where $t$ is the approximation rate for functions in $W$, as in \eqref{ApproxProp}.
\end{theorem}

\beginproof
Similar to the proof of Theorem \ref{SDError}, we write
\begin{equation*}
\|u(y)-u_{h}(y)\|^2_V\leq \exp\left(2\|b(y)\|_{L^\infty}\right)\underset{v_h\in V_h}{\min}\|u(y)-v_{h}\|^2_{V}.
\end{equation*}
Using the approximation property \eqref{ApproxPropalpha} applied to $u(y)$ we obtain
\begin{equation*}
\|u(y)-u_{h}(y)\|^2_V\leq C\exp\left(2\|b(y)\|_{L^\infty}\right)n_h^{-2 t \alpha}\|u(y)\|^2_{W^\alpha}, 
\end{equation*}
and by Cauchy-Schwarz inequality
\begin{align}
\|u-u_{h}\|^2_{\cV_2}\leq Cn_h^{-2 t \alpha}\mathbb{E}(\exp(4\|b(y)\|_{L^\infty})^{\frac{1}{2}} \|u\|^2_{\cW_4^\alpha}.
\end{align}
This concludes the proof of the theorem.
\endproof


We close this subsection by discussing conditions that ensure that 
the Gaussian field $b$ belongs to the H\"older space $C^{\alpha}(D)$ almost surely,
with bounded moments 
\be
\mathbb{E}(\|b\|^k_{C^{\alpha}})
<\infty, \quad k<\infty,
\label{moments}
\ee
as needed in Proposition \ref{spaceregrough}. 

We first observe that the boundedness of all moments holds provided that the Gaussian field $b$ 
belongs to $C^{\alpha+\e}(D)$ almost surely for some $\e>0$. Indeed, this implies that $b$ 
belongs to the space $C_0^{\alpha}(D)$, which is defined by the condition
\begin{equation*} 
\underset{x_1\rightarrow x_2}{\lim}\frac{|b(x_1)-b(x_2)|}{|x_1-x_2|^{\alpha}}=0.   
\end{equation*}
This space is a separable Banach space, endowed with the same norm as $C^\alpha(D)$. Therefore Fernique's theorem 
implies that the moment bound  \eqref{moments} holds.
 
There are several criteria which ensure that $b$ belongs to the H\"older space $C^{\alpha}(D)$. A standard probabilistic one is Kolmogorov's continuity theorem (see e.g. Theorem $2.1$ \cite{RY13}),
that reads on the covariance $C_b(x,x'):=\E(b(x)b(x'))$ of the Gaussian random field: if for some $K>0$ and $\beta >\alpha$,
\begin{align}
C_b(x,x)+C_b(x',x')-2C_b(x,x')\leq K|x-x'|^{2\beta}, \quad x,x'\in D,
\label{covar}
\end{align}
then $b$ admits a continuous modification which belongs to $C^{\alpha}(D)$ almost surely.
In particular, this also implies that $b$ belongs to $C^{\alpha+\e}$ for some $\e>0$ sufficiently small 
which thus implies the moment bounds \eqref{moments}. Note that \eqref{covar} holds in particular when the covariance
function $C_b$ has
H\"older smoothness $C^{2\beta}(D\times D)$. We apply this criterion in \S 6 to treat the case of
Brownian and Mat\'ern processes.

\section{Weighted least-squares method}
\label{sec:wls}
\noindent 
For a given set $\Lambda_n\subset \cF$ with $n=\#(\Lambda_n)$,
the weighed least-squares methods allows us to construct an approximation
of the solution map $u$ in $\cV_{\Lambda_n}$ 
from a finite number of evaluations of $u$: we are given $m$ evaluations 
$u^1,\ldots,u^m$ of $u$ at the points $y^1,\ldots,y^m \in U$,
that is,
\be
\label{eq:obs_model}
u^i=u(y^i), \quad i=1,\dots,m,
\ee
and define the weighted least-squares estimator $u_n^W\in \cV_{\Lambda_n}$ of $u$ as 
\begin{equation}
u_n^W:= \argmin_{v \in \cV_{\Lambda_n}}  \sum_{i=1}^m w^i \| v(y^i) - u^i \|^2_V,
\label{eq:def_wls}
\end{equation}
where $w^i\geq 0$ are positive weights. The properties of stability and accuracy of this estimator 
depend on the distribution of the samples $y^1,\ldots,y^m$, on their number $m$ and on the 
choice of the weights. 

\subsection{Optimal weights and sampling measure}

Following the analysis of weighted least squares in \cite{CM2016}, we introduce the  nonnegative weight function 
\begin{equation}
y\mapsto w(y):=
\dfrac{n}{
\sum_{\nu \in \Lambda_n} 
| H_{\nu}(y)|^2
},
\label{eq:def_weight} 
\end{equation}
for a given $\Lambda_n$ associated to $V_n$, where we recall that $H_\nu$ are orthonormal in $L^2(U,d\gamma)$.  
Using the function \eqref{eq:def_weight} we define the probability measure 
\begin{equation}
d\mu := w^{-1} d\gamma = 
\dfrac{1}{n}
\sum_{\nu \in \Lambda_n} 
| H_{\nu}(y)|^2 d\gamma. 
\label{eq:def_sampling_measure}
\end{equation}
The $y^1,\ldots,y^m$ used in \eqref{eq:def_wls} are independent samples drawn from $d\mu$.  
The evaluations of the function \eqref{eq:def_weight} at the selected points 
$$w^i:=w(y^i),$$
give the weights in \eqref{eq:def_wls}. 
By expanding the estimator $u_n^W$ over the orthonormal basis we obtain  
$$
u_n^W=\sum_{\nu \in \Lambda_n} v_\nu H_\nu, 
$$
where $\bv:=(v_\nu)_\nu$ denotes the collection of the coefficients. 
The calculation of the estimator $u_W$ is equivalent to solving the normal equations 
\begin{equation}
\bG \bv = \bd,
\label{eq:sys}
\end{equation}
where $\bG$ and $\bd$ are defined element-wise as 
\begin{equation}
\bG_{\nu,\widetilde{\nu}}:= \dfrac{1}{m} \sum_{i=1}^m w^i H_\nu(y^i)  H_{\widetilde{\nu}}(y^i),
\qquad 
\bd_{\nu}:= \dfrac{1}{m} \sum_{i=1}^m w^i u^i H_\nu(y^i). 
\label{eq:def_G_d}
\end{equation}

The solution to the weighted least-squares always exists and is unique if and only if 
$\bG$ is non-singular. In order to avoid the situation where $\bG$ is ill-conditionned,
we consider the conditioned least-squares estimator $u_n^C$ 
introduced in \cite{CM2016}, that is defined as 
$$
u_n^C:= 
\begin{cases}
u_n^W, & \textrm{ if } \| \bG- \bI \|_2 \leq \frac12, \\
0, & \textrm{ otherwise}.
\end{cases}  
$$
The choice of the particular weight function \eqref{eq:def_weight} and sampling measure
\eqref{eq:def_sampling_measure} is crucial, in that it yields stability and optimal accuracy
of the method when the number $m$ of sample to be evaluated is larger than $n$ only by a logarithmic factor.
This is expressed by the following result, which is an adaptation of Theorem 13 in \cite{CM2016b} 
to the Hermite polynomial case.

\begin{theorem}
\label{thm:wls_estimator}
For any real $s>0$ and any integer $n\geq 1$, if $m$ satisfies  
\begin{equation}
n\leq \kappa \dfrac{m}{\ln m}, \textrm{ with } \kappa=\kappa(s):=\dfrac{1-\ln 2}{2 + 4s},
\label{eq:condmn_wls}
\end{equation}
and the samples $y^1,\ldots,y^m$ are independent and 
drawn from \eqref{eq:def_sampling_measure} 
then the estimator $u_n^C$ satisfies 
\begin{equation}
\mathbb{E}(\|  u - u_n^C \|^2_{\cV_2} ) \leq 
\left( 1+\beta(m)  \right) 
\min_{ v \in \cV_{\Lambda_n} } \| u - v \|^2_{\cV_2}
+ 2 \| u \|^2_{\cV_2} m^{-2s}, 
\label{eq:estimate_wls}
\end{equation}
with $\beta(m):= \dfrac{4\kappa}{\ln m}\to 0$ as $m\to +\infty$. 
\end{theorem}

This allows us to obtain convergence estimates in expectation for the solution to \eqref{eq:model_pde}.

\begin{theorem}
\label{thm:wls_for_lognormal_pb}
Let $(\rho_j)_{j\geq 1}$ be a sequence that satisfies the assumptions of Theorem \ref{thm:theoherm},
and let $\Lambda_n$ be sets corresponding to the $n$ largest 
$\xi_\nu^{-1/2}$, with $\xi_\nu$ given by~\eqref{eq:weights_lognormal}. Then the estimator $u_n^C$ 
built with $m$ samples drawn from \eqref{eq:def_sampling_measure}
 under  \eqref{eq:condmn_wls} with $s:=1/q$
satisfies 
\begin{equation}
\mathbb{E}( \| u-u_n^C \|^2_{\cV_2} ) \leq 
\left( 
C \left( 
1+\dfrac{ 4 \kappa(2s) }{\ln n}
\right)
+2 \| u \|^2_{\cV_2}
\right) 
n^{-2s},
\label{eq:est_wls_lognormal}
\end{equation}
where the constant $C$ is the same as in Theorem~\ref{thm:theoherm}. 
\end{theorem} 

The convergence estimate \eqref{eq:est_wls_lognormal}
 shows that, in the ideal situation where exact evaluations of $u$ are available, 
%
the approach based on weighted least squares achieves in expectation the optimal convergence rate $n^{-s}$
for the  approximation error in $\cV_2$ of the solution to the  lognormal PDE of interest.

In practice, the exact evaluation $u^i=u(y^i)$ is not available and is replaced 
by its finite element approximation $u^i_h:=u_h(y^i)$. 
Following the approach outlined in \S 2, we 
thus apply the weighted least-squares approximation 
to the map $y\mapsto u_h(y)$, and denote by $u^C_{h,n}\in V_h\otimes \P_{\Lambda_n}$
the resulting estimator of $u_h$. By the same argument, we reach the following analog of Theorem \ref{thm:wls_for_lognormal_pb}.

\begin{theorem}
\label{thm:wls_for_lognormal_pbh}
Let $(\rho_j)_{j\geq 1}$ be a sequence that satisfies the assumptions of Theorem \ref{thm:theoherm},
and let $\Lambda_n$ be sets corresponding to the $n$ largest 
$\xi_\nu^{-1/2}$, with $\xi_\nu$ given by~\eqref{eq:weights_lognormal}. Then the estimator $u_{h,n}^C$ 
built with $m$ samples $u^i_h:=u_h(y^i)$ drawn from \eqref{eq:def_sampling_measure}
 under  \eqref{eq:condmn_wls} with $s:=1/q$
satisfies 
\begin{equation}
\mathbb{E}( \| u_h-u_{h,n}^C \|^2_{\cV_2} ) \leq 
\left( 
C \left( 
1+\dfrac{ 4 \kappa(2s) }{\ln n}
\right)
+2 \| u_h \|^2_{\cV_2}
\right) 
n^{-2s},
\label{eq:est_wls_lognormalh}
\end{equation}
where the constant $C$ is the same as in Theorem~\ref{thm:theoherm}. 
\end{theorem} 

By Cauchy-Schwarz inequality, the estimate in the above theorem implies
\begin{equation*}
\mathbb{E}( \| u_h-u_{h,n}^C \|_{\cV_2} ) \leq Cn^{-s}, 
\end{equation*}
with $s:=1/q$.
Combining this result together with Theorem \ref{SDError} or with Theorem \ref{thm:RoughError}, we obtain an estimate on the total error,
\be
\label{eq:toterrwls}
\E(\|u-u_{h,n}^C\|_{\cV_2})\leq C(n_h^{-\alpha t}+n^{-s})
\ee
where $C>0$ and $\alpha$ is equal to $1$ if the Gaussian field $b$ is smooth (Theorem \ref{SDError}) 
or is equal to the global H\"older regularity index of b if it is rough (Theorem \ref{thm:RoughError}). 


\subsection{Truncation to a finite number of random variables}\label{sec:truncFiniteNumberRV}
\noindent
At this point, it is worth to observe that $u_h(y^i)$ 
still cannot be evaluated,  
due to the infinite number of coordinates of $y^i$ that prevents the practical generation of the samples from \eqref{eq:def_sampling_measure}. 
This difficulty can be overcome by using a proper truncation of 
the coordinates of $y^i$.


For some finite $J$ to be fixed further, we truncate $y^i$ by setting to $0$ 
its coordinates $y^i_j$ for all $j>J$. The resulting samples $y^{1,J},\ldots,y^{m,J}$ can be directly generated from the probability measure 
\begin{equation}
d\mu_J :=  
\dfrac{1}{n}
\(\sum_{\nu \in \Lambda_n} 
\prod_{1\leq j \leq J} | H_{\nu_j}(y_j)|^2\)d\gamma_J, 
\label{eq:def_sampling_measure_truncated}
\end{equation}
where $\gamma_J$ is the $J$-dimensional Gaussian measure. 
We thus now evaluate 
\begin{equation*}
u_{h,J}^i=u_h(y_J^i).
\end{equation*}
The application of the weighted-least square method to these data yields an approximation
$u^C_{h,n,J}$ to the function $u_{h,J}$ defined in a similar way as $u_h$ by \eqref{SpaceDiscrPara}, however
with diffusion coefficient $a$ replaced by
\begin{align}
\label{eq:b_J}
a_J(y):=\exp(b_J(y)), \quad b_J(y):=\sum_{j=1}^{J}y_j\psi_j(x).
\end{align}

The following result is obtained by the exact same arguments as those leading to Theorem~\ref{thm:wls_for_lognormal_pbh}.

\begin{theorem}
\label{finalthm_wls}
Let $(\rho_j)_{j\geq 1}$ be a sequence that satisfies the assumptions of Theorem \ref{thm:theoherm},
and let $\Lambda_n$ be sets corresponding to the $n$ largest 
$\xi_\nu^{-1/2}$, with $\xi_\nu$ given by~\eqref{eq:weights_lognormal}. 
Then the estimator $u^C_{h,n,J}$ built with $m$ samples 
$y_J^1,\ldots,y_J^m$ 
drawn from \eqref{eq:def_sampling_measure_truncated}
 under  \eqref{eq:condmn_wls} and with the observation model $u_{h,J}^i=u_h(y_J^i)$ 
satisfies 
\begin{align*}
\mathbb{E}( \| u_{h,J}-u^C_{h,n,J} \|^2_{\cV_2} ) \leq 
\left( 
C \left( 
1+\dfrac{ 4 \kappa(2s) }{\ln n}
\right)
+2 \| u \|_{\cV_2}^2
\right) 
n^{-2s},
\end{align*}
where $s=1/q$ and the constant $C$ is the same as in Theorem~\ref{thm:theoherm}.
\end{theorem}
\noindent
In order to estimate the total error $\| u - u^C_{h,n,J} \|_{\cV_2}$, we write
\begin{equation*}
\| u - u^C_{h,n,J} \|_{\cV_2} \leq 
\| u - u_h \|_{\cV_2}
+
\| u_h - u_{h,J} \|_{\cV_2}
+
\| u_{h,J} - u_{C} \|_{\cV_2}.
\end{equation*}
The first term is estimated by Theorem \ref{SDError} or Theorem \ref{thm:RoughError}, and the last one
by Theorem \ref{finalthm_wls}. We are thus left with estimating the second term $\| u_h - u_{h,J} \|_{\cV_2}$
that describes the error when the variables $y_j$ are put to $0$ for $j>J$. This is the object of the following result.

\begin{prop}
\label{eq:truncation_error_wls}
Let $(\rho_j)_{j\geq 1}$ be a sequence of positive real numbers such that \eqref{rhocond} holds and that $\sum_{j\geq 1}\exp(-\rho_j^2)<+\infty$ . Assume in addition that $\(\frac {\sqrt{\log(1+j)}}{\rho_j}\)_{j\geq 1}$ is a 
non-increasing sequence. Then, for any $J\geq 1$,
\begin{align}
\label{uhuhJ2}
\|u_h-u_{h,J}\|_{\cV_2} \leq &
C_3\frac{\sqrt{\log(1+J)}}{\rho_J},
\end{align}
with $C_3>0$.
\end{prop}





\beginproof
Since $u_h$ and $u_{h,J}$ are the solutions to \eqref{SpaceDiscrPara} with diffusion coefficients $a=\exp(b)$ and $a_J=\exp (b_J)$, respectively, we derive by substraction of the two variational formulations the standard stability estimate
\begin{align*}
\|u_h(y)-u_{h,J}(y)\|_{V}&\leq \exp\left(\|b_J(y)\|_{L^\infty}\right)\|a(y)-a_J(y)\|_{L^\infty}\|u_h(y)\|_V.
\end{align*}
Using $|e^x-e^y|\leq |x-y|(e^x+e^y)$ for any $x,y\in\mathbb{R}$, 
\begin{align*}
\|u_h(y)-u_{h,J}(y)\|_{V}&\leq \exp\left(\|b_J(y)\|_{L^\infty}\right) \|b(y)-b_J(y)\|_{L^\infty}\|u_h(y)\|_V 
\(\exp(\|b(y)\|_{L^\infty})+\exp(\|b_J(y)\|_{L^\infty})\).
\end{align*}
From \eqref{rhocond}, it is readily checked that $\|b(y)-b_J(y)\|_{L^\infty}\rightarrow 0$ when $J\rightarrow \infty$ for any $y$ such that $\underset{j\geq 1}{\sup}\frac{|y_j|}{\rho_j}<\infty$. Then, for such $y$, using the 
a-priori estimate for $\|u_h(y)\|_V$, we obtain
\begin{align}\label{DiffL}
\|u_h(y)-u_{h,J}(y)\|_{V}\leq 2\|f\|_{V'}\exp\left(3K\ \underset{j\geq 1}{\sup} \frac{|y_j|}{\rho_j} \right)\|b(y)-b_J(y)\|_{L^\infty}\end{align}
where $K$ is as in \eqref{rhocond}. We bound the term $\|b(y)-b_J(y)\|_{L^\infty}$ by writing, for all $x\in D$,
\begin{align}\label{DiffField}
|b(x,y)-b_J(x,y)| &\leq \underset{j\geq 1}{\sup}\left|\frac{y_j}{\sqrt{\log(1+j)}}\right|\sum_{j\geq J+1}\sqrt{\log(1+j)} |\psi_j(x)| \nonumber\\
&\leq K \frac{\sqrt{\log(1+J)}}{\rho_J} \underset{j\geq 1}{\sup}\left| \frac{y_j}{\sqrt{\log(1+j)}} \right|,
\end{align}
where the second inequality uses the assumption that $\(\frac {\sqrt{\log(1+j)}}{\rho_j}\)_{j\geq 1}$ is 
non-increasing.

In order to pursue, we prove that $\underset{j\geq 1}{\sup}\left| \frac{y_j}{\sqrt{\log(1+j)}} \right|$ is bounded by a quantity $C(y)$ which admits moments of any order with respect to the Gaussian measure $d\gamma$. We follow closely the proof of Lemma $1$ in \cite{AT03}. With $\epsilon_j:=\sqrt{\log(1+j)}$, for all $j\geq 1$, let $b>0$ be a deterministic constant to be chosen later. By standard estimates,
\begin{align*}
\int_{|y_j|\geq b\epsilon_j} d\gamma(y)\leq \frac{C}{(j+1)^{\frac{b^2}{2}}}.
\end{align*}
Choosing $b>\sqrt{2}$, Borel-Cantelli lemma ensures the existence of a subset $\tilde{U}\subset U$ with full measure 
such that for all $y\in \tilde{U}$, there exists $j_0(y)\geq 1$ such that
\begin{align*}
|y_j|\leq b\sqrt{\log(1+j)}, \quad j\geq j_0(y).
\end{align*}
Now, we may define the random variables
\begin{equation*}
\tau(y):=\min \{j\geq 1:\ |y_j|\leq b\sqrt{\log(1+j)}\},
\end{equation*}
and
\begin{equation*}
D(y):=\underset{1\leq j\leq \tau(y)}{\sup}|y_j|.
\end{equation*}
Thus, for $j\geq 1$
\begin{align}\label{ineq:sup}
|y_j| \leq C(y)\sqrt{\log(1+j)},\quad\quad y\in \tilde{U},
\end{align}
with $C(y)=b\log(2)\max(D(y),1)$. Combining \eqref{DiffL}, \eqref{DiffField} and \eqref{ineq:sup}, for $y\in \tilde{U}$
\be
\label{uhuhJ}
\|u_h(y)-u_{h,J}(y)\|_{V}\leq 2K\|f\|_{V'} \frac{\sqrt{\log(1+J)}}{\rho_J} C(y) \exp\left(3K\underset{j\geq 1}{\sup}\left|\frac{y_j}{\rho_j}\right|\right).
\ee
From the results in \cite{BCDM}, we know that \eqref{rhocond} with $\sum_{j\geq 1}\exp(-\rho_j^2)<+\infty$ 
implies that 
\begin{equation*}
\int_U\exp\left(p\; \underset{j\geq 1}{\sup}\left|\frac{y_j}{\rho_j}\right|\right)d\gamma(y)<\infty,
\end{equation*}
for all $p>0$. In addition, the proof of Lemma $1$ in \cite{AT03}, shows that for $p\geq 1$ and $b>2\sqrt{2}$,
\begin{align*}
\int_{\tilde{U}}|D(y)|^p d\gamma(y)<\infty.
\end{align*}
Therefore \eqref{uhuhJ2} follows by integration of \eqref{uhuhJ} after squaring, and by application of Cauchy-Schwarz inequality.
\endproof
\noindent

Let us observe that by a proper reordering of the variables, we may always assume that $\(\frac {\sqrt{\log(1+j)}}{\rho_j}\)_{j\geq 1}$ is 
non-increasing, so that the above proposition can be applied. In addition, if $(\rho_j^{-1})_{j \geq 1} \in \ell^q(\mathbb{N})$, we may write
\begin{equation*}
\(\frac {\sqrt{\log(1+J)}}{\rho_J}\)^q\leq J^{-1} 
\sum_{j=1}^J\(\frac {\sqrt{\log(1+j)}}{\rho_j}\)^q \leq 
J^{-1}(\log(1+J))^q \|(\rho_j^{-1})_{j \geq 1} \|_{\ell^q}^q.
\end{equation*}

Combining the previous results, we conclude this section with the following error estimate for the weighted least-squares
method with space discretization and truncation of the variables.
\begin{theorem}
\label{finalthm_wls_inf_dim}
Let $(\rho_j)_{j\geq 1}$ be a sequence that satisfies the assumptions of Theorem \ref{thm:theoherm},
and let $\Lambda_n$ be sets corresponding to the $n$ largest 
$\xi_\nu^{-1/2}$, with $\xi_\nu$ given by~\eqref{eq:weights_lognormal}. Then the estimator $u^C_{h,n,J}$ built with $m$ samples 
$y_J^1,\ldots,y_J^m$ 
drawn from \eqref{eq:def_sampling_measure_truncated}
 under  \eqref{eq:condmn_wls} and with the observation model $u_{h,J}^i=u_h(y_J^i)$ 
satisfies 
\begin{align*}
\mathbb{E}( \| u-u^C_{h,n,J} \|_{\cV_2} ) \leq &
C\(n_h^{-\alpha t} + n^{-s}+(\log(1+J))^2J^{-s}\),
\end{align*}
where $C>0$, $s=1/q$, and where $\alpha$ is equal to $1$ if the Gaussian field $b$ is smooth (Theorem \ref{SDError}) 
or is equal to the global H\"older regularity index of b if it is rough (Theorem \ref{thm:RoughError}). 
\end{theorem}

From 
the above estimate, 
a reasonable strategy is 
to choose the truncation level $J$ of the same order as $n$, up to logarithmic factors.


\section{Examples of random fields}
\label{ref:examples_random_field}
The analysis in \S\ref{sec:SD} and  
\S\ref{sec:wls} shows that
for both Galerkin and weighted least-squares methods, the total
error of approximation is split into two terms, resulting from
the spatial and parametric discretization, respectively. 
The spatial error term is controlled by the H\"olderian regularity of the sample path of $b$,
while the parametric error term is controlled by the size properties
of the functions $\psi_j$ in the representation of the Gaussian random field $b$. 

In the following subsections we compute more explicitly the exponents appearing in 
\eqref{eq:toterrwls} 
for relevant Gaussian field: Brownian bridge on $[0,1]$ and Gaussian fields
with Mat\'ern covariances. As to the functions $(\psi_j)_{j\geq 1}$
we consider wavelet representations adapted to such fields, as studied in \cite{BCM}.
 We show that, for these relevant 
examples, the spatial discretization and parametric rate exponents coincide. This leads
us to error estimates in terms of the total number of degrees of freedom
\begin{equation*}
n_{dof}:=n_h n,
\end{equation*}
that describes the approximant $u_{h,n}^C$. 
Here, we do not discuss anymore the third error term due to variable truncation that
appears in Theorem \ref{finalthm_wls_inf_dim}, which is specific to the
weighted least-squares method, and can always be absorbed in the second term by taking $J$ slightly greater
than $n$.

\subsection{Brownian bridge}
In this subsection, we consider Brownian type processes. For simplicity, we focus on Brownian bridge but a similar analysis can be performed for Brownian motion. First, recall that a Brownian bridge on $[0,1]$ is the centered Gaussian process starting from $0$ with covariance function given by
\begin{equation*}
C_b(x_1,x_2):=x_1\wedge x_2-x_1x_2,\quad x_1,x_2\in [0,1].
\end{equation*}
It is well known that, by application of Kolmogorov's continuity theorem, for any $\delta<1/2$, there exists a continuous modification such that its sample paths have $C^{\delta}$ H\"older regularity on $[0,1]$ almost surely. In addition, one has
$\E(\|b\|_{C^\delta}^k)<\infty$ for all $k<\infty$. Then, Theorem \ref{thm:RoughError} applies with $\alpha=1/2-\varepsilon$ for any $0<\varepsilon<1/2$, and $t=1$ since 
we work here in spatial dimension $d=1$. This leads to the spatial discretization error
\begin{equation*}
\|u-u_{h}\|_{\cV_2}\leq Cn_h^{-(\frac{1}{2}-\epsilon)},\quad\quad h>0.
\end{equation*}
It is well known that Brownian bridge admits several series expansion. As discussed in \cite{BCDM},
the most relevant one for the purpose of polynomial approximation is in terms of the Schauder basis
\be
b(x):=\sum_{l\geq 0}\sum_{k=0}^{2^l-1}y_{l,k}\psi_{l,k}(x),
\label{eq:schauder_basis}
\ee
where the $y_{l,k}$ are independent standard normal random variables.
The functions $\psi_{l,k}$ are defined by
\be
\label{eq:def_schauder_basis_bb}
\psi_{l,k}(x):=2^{-\frac{l}{2}}\psi(2^l x-k),\quad x\in [0,1]
\ee
with $\psi(x)=\max\{1/2-|x-1/2|,0 \}$. This is  also known as the Levy-Ciesieslki
representation, and can be rewritten as $\sum_{j\geq 1} y_j\psi_j$, by
enumeration of all functions from coarser to finer scales, \emph{i.e.}~$j=2^l+k$.

Concerning the parametric discretization error, as observed in \cite{BCDM}, from the decay and localization properties
of the Schauder basis functions, one has
\begin{equation*}
\sum_{l\geq 0}\sum_{k=0}^{2^l-1} 2^{l \beta} |\psi_{l,k}(x)|<\infty,\quad x\in [0,1]
\end{equation*}
for all $\beta<1/2$, or equivalently
\be
\label{eq:weights_brownian_bridge}
\sum_{j\geq 1} j^\beta |\psi_{j}(x)|<\infty, \quad x\in [0,1].
\ee
This shows that the assumptions of Theorem \ref{thm:theoherm} are satisfied 
for all $q>1/\beta>2$, that is, $q$ can be made arbitrarily close to $2$.
Thus for any $0< \e <1/2$, 
Theorem~\ref{thm:wls_for_lognormal_pbh}
gives the estimate
\be
\E(\|u_h-u_{h,n}^C\|_{\cV_2})\leq C n^{-(\frac{1}{2}-\e)},
\label{eq:est_wls_bb}
\ee
for the weighted least-squares approximation. 
The rate of decay in terms of $n_h$ and $n$ coincide for the spatial and parametric error terms.
This suggests to take $n=n_h$, which leads, for any $\e>0$ to 
a global error estimate $n_{dof}^{-(\frac{1}{4}-\varepsilon)}$ in terms of the total number
of degrees of freedom.

\subsection{Stationary Gaussian Mat\'ern fields}

In this subsection, we consider the centered Gaussian Mat\'ern fields on a general
domain $D\subset \R^d$. The covariance function for such fields is given by
\begin{equation*}
\mathbb{E}(b(x_1)b(x_2))=C_{\delta,\ell}(|x_1-x_2|), \quad x_1,x_2\in D,
\end{equation*}
with 
\begin{equation*}
C_{\delta,\ell}(r)=\dfrac{2^{1-\delta}}{\Gamma(\delta)}\left(\dfrac{\sqrt{2\delta}r}{\ell}\right)^{\delta}K_\delta\left(\dfrac{\sqrt{2\delta}r}{\ell}\right)
\end{equation*}
where $\ell,\delta>0$, $K_{\delta}$ is the modified Bessel function of the second kind and $\Gamma$ is the gamma function. We restrict our analysis to small parameters values $\delta< 1$, that correspond to rough fields, but a similar one can be performed for $\delta\geq 1$.  Based on series expansion of the modified Bessel function $K$ when $r\rightarrow 0^+$, 
\begin{equation*}
|1-C_{\delta,\ell}(r)|\leq C r^{2\delta}
\end{equation*}
for some positive constant $C$ only depending on $\delta$ and $\ell$. Then, for all $(x_1,x_2)\in D\times D$,
\begin{equation*}
2(1-\mathbb{E}b(x_1)b(x_2))\leq C|x_1-x_2|^{2\delta}
\end{equation*}
for some $C>0$ only depending on $\delta$ and $\ell$. Thus, Kolmogorov's continuity theorem implies that the centered Mat\'ern Gaussian random fields admit continuous modifications which belong to $C^{\delta'}(\overline{D})$ almost surely for any $\delta'<\delta$. Then, Theorem \ref{thm:RoughError} applies with $\alpha=\delta-\varepsilon$ for any $0<\varepsilon<\delta$. Moreover, when $W$ is equal to $V\cap H^2(D)$,  $t$ is equal to $1/d$ leading to the following spatial discretization error,
\begin{equation*}
\|u-u_{h}\|_{\cV_2}\leq Cn_h^{-(\frac{\delta}{d}-\varepsilon)},\quad\quad h>0.
\end{equation*}
Based on wavelet expansion obtained in \cite{BCM} for the Gaussian Mat\'ern fields, it is possible to compute explicitly a parametric discretization rate. For this purpose, let $(\psi_\lambda)_{\lambda\in \mathcal{I}}$ be the wavelet basis defined in Section $4$ of \cite{BCM} where $\lambda$ is a scale-space index with $|\lambda|$ denoting the scale parameter and where $\mathcal{I}$ denotes the set of these indices with $\# \{\lambda\in\mathcal{I}:\, |\lambda|=l \}\sim 2^{dl}$ for $l\geq 0$. The following crucial localization property is ensured by Corollary $4.3$ of \cite{BCM}
\begin{equation*}
\underset{x\in D}{\sup} \sum_{|\lambda|=l}|\psi_\lambda(x)|\leq C 2^{-\delta l},\quad l\geq 0
\end{equation*}
for some constant $C>0$ independent of $l$. Then Corollary $1.3$ in \cite{BCM} can be applied and  parametric rates $s$ strictly less than $\delta/d$ are achievable. 
Setting $\rho_\lambda=2^{|\lambda| \epsilon}$ 
with $\lambda\in\mathcal{I}$, 
for any $0< \e <\delta/d$, Theorem \ref{thm:wls_for_lognormal_pbh}
gives the estimate
\begin{equation*}
\E(\|u_h-u_{h,n}^C\|_{\cV_2})\leq C n^{-(\frac{\delta}{d}-\e)},
\end{equation*}
for the weighted least-squares approximation. We note
that the rate of decay in terms of $n_h$ and $n$ coincide for the spatial and parametric error terms.
This suggests to take $n=n_h$, which leads for any $\e>0$ to 
a global error estimate $n_{dof}^{-(\frac{\delta}{2d}-\varepsilon)}$ in terms of the total number of degrees of freedom.


\section{Numerical illustration with the Brownian bridge}
\label{sec:num_illust}
In this section we numerically verify the error estimate 
\eqref{eq:est_wls_bb} for the weighted least-squares method.  
%
The estimator $u_{h,n}^C$ is calculated using $m=3 \, n \, \lceil \ln n \rceil$ 
random samples. 
The index set $\Lambda_n$ is chosen by taking the $n$ smallest $\xi_\nu$ in \eqref{eq:weights_lognormal} with $r=1$, 
as an application of 
Theorems~\ref{summability_theo_lognormal}, \ref{thm:lemmaqmh}, \ref{thm:theoherm}. 
The weights $\xi_\nu$ are calculated 
by \eqref{eq:weights_brownian_bridge} 
from the sequence $(\rho_i)_{i\geq 1}$ 
defined as 
$$
\rho_i := \tilde{\rho}_i \dfrac{ \ln 2}{2C\sqrt{r}}, \qquad 
\tilde\rho_i:= 2^{\beta  \lfloor \log_2 i \rfloor },  \qquad 
C:= \sup_{x \in D} \sum_{i=1}^{2^{L+1}-1} \tilde{\rho}_i  | \psi_i(x) |   
$$ 
 for any $\beta \in (0,\frac12)$, that satisfies 
\eqref{eq:psicr} as a finite sum up to $2^{L+1}-1$ and with an additional factor $\frac12$. 
The set $\Lambda_n$ produced with this construction therefore depends on the value of $\beta $: large values of $\beta$ promote higher degrees along the first levels of the Schauder basis, 
and small values of $\beta$ generate more isotropic spaces. 

If the Schauder basis $\psi_i$ is rescaled as $\tau \psi_i$ then \eqref{eq:psicr} still holds with the rescaled sequence of weights $(\tau^{-1}\rho_i)_i$.  
Rescaling with $\tau<1$ reduces the variance of the random field 
\eqref{eq:schauder_basis}, 
and also promotes higher degrees along the first levels of the Schauder basis. 
It is worth to notice that 
the convergence estimate \eqref{eq:est_wls_bb}
requires $\beta <\frac12$,    
but 
similar polynomial spaces as those generated with any value of $\beta$  
(even larger than $\frac12$) 
can actually be generated 
taking a smaller value of $\beta$ and rescaling the Schauder basis by a constant.
Therefore, in the following we only verify that 
the numerical convergence rates agree with \eqref{eq:est_wls_bb}, 
but we do not advocate a precise value of $\beta$ to use in the computations, for example the one giving the fastest convergence rate.
We test $\tau=1$ and 
$\beta=\frac18 $, 
$\beta=\frac14$, $\beta=\frac12$.

\begin{figure}[!h]
\center
\includegraphics[scale=0.4]{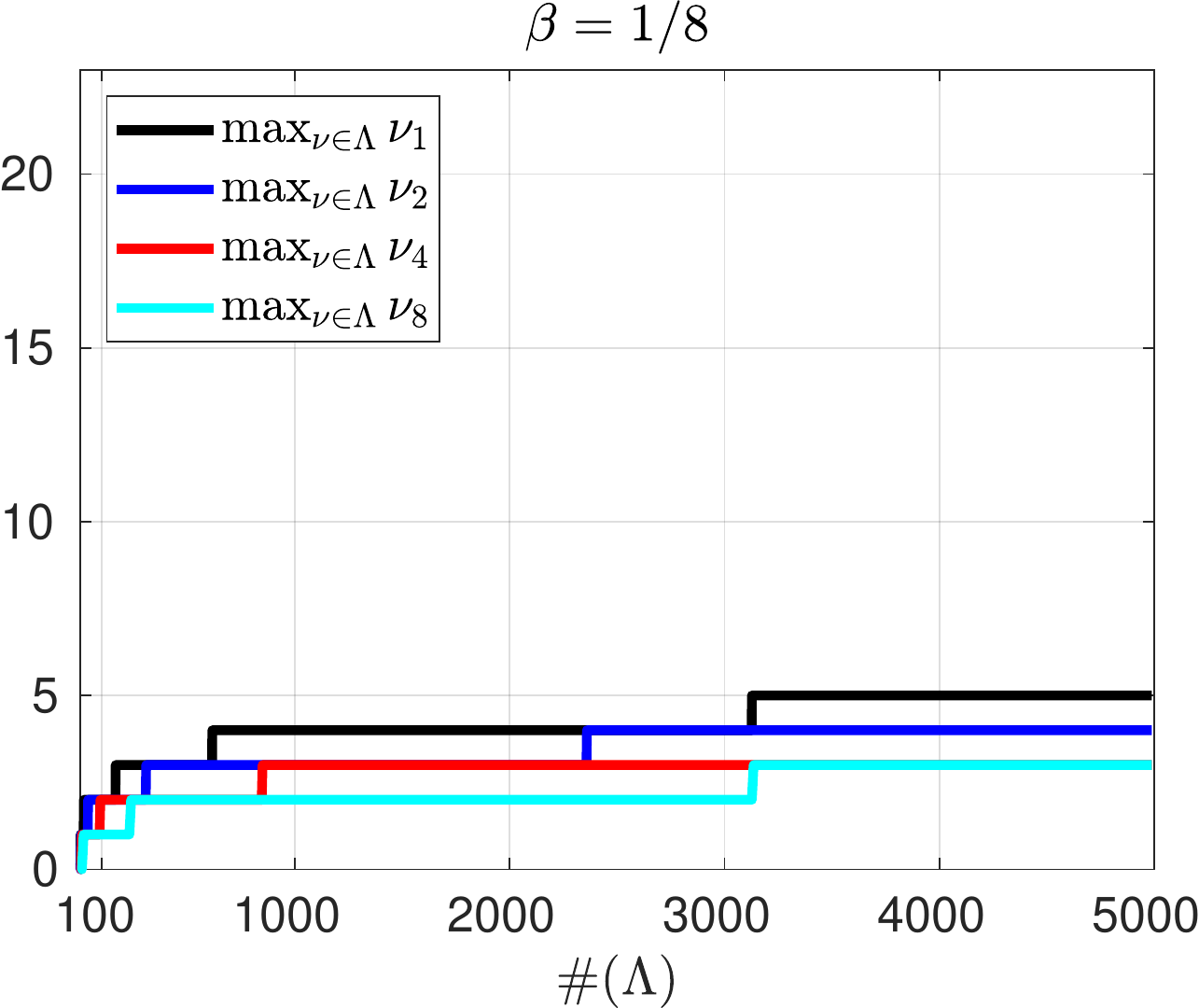} 
\includegraphics[scale=0.4]{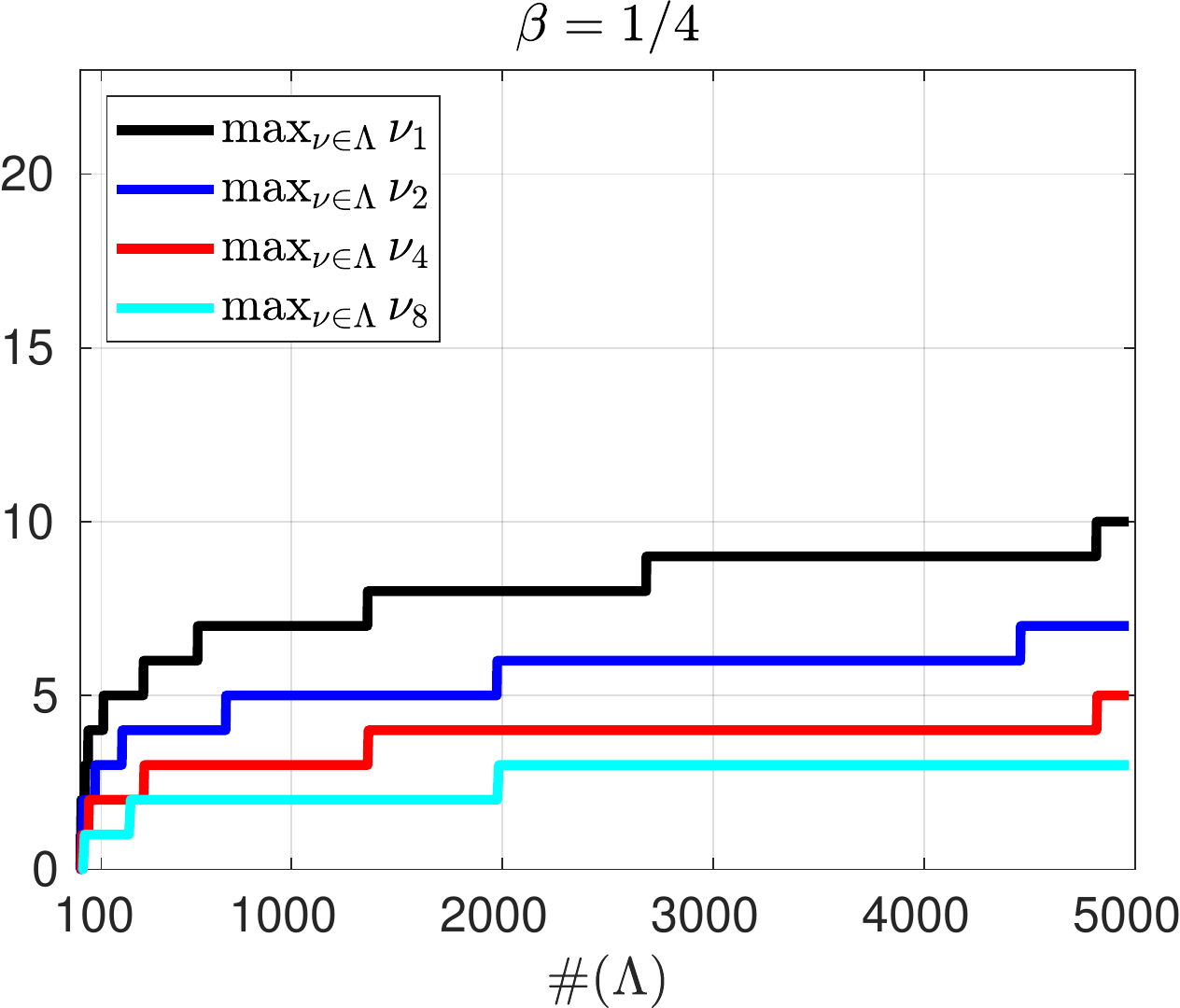} 
\includegraphics[scale=0.4]{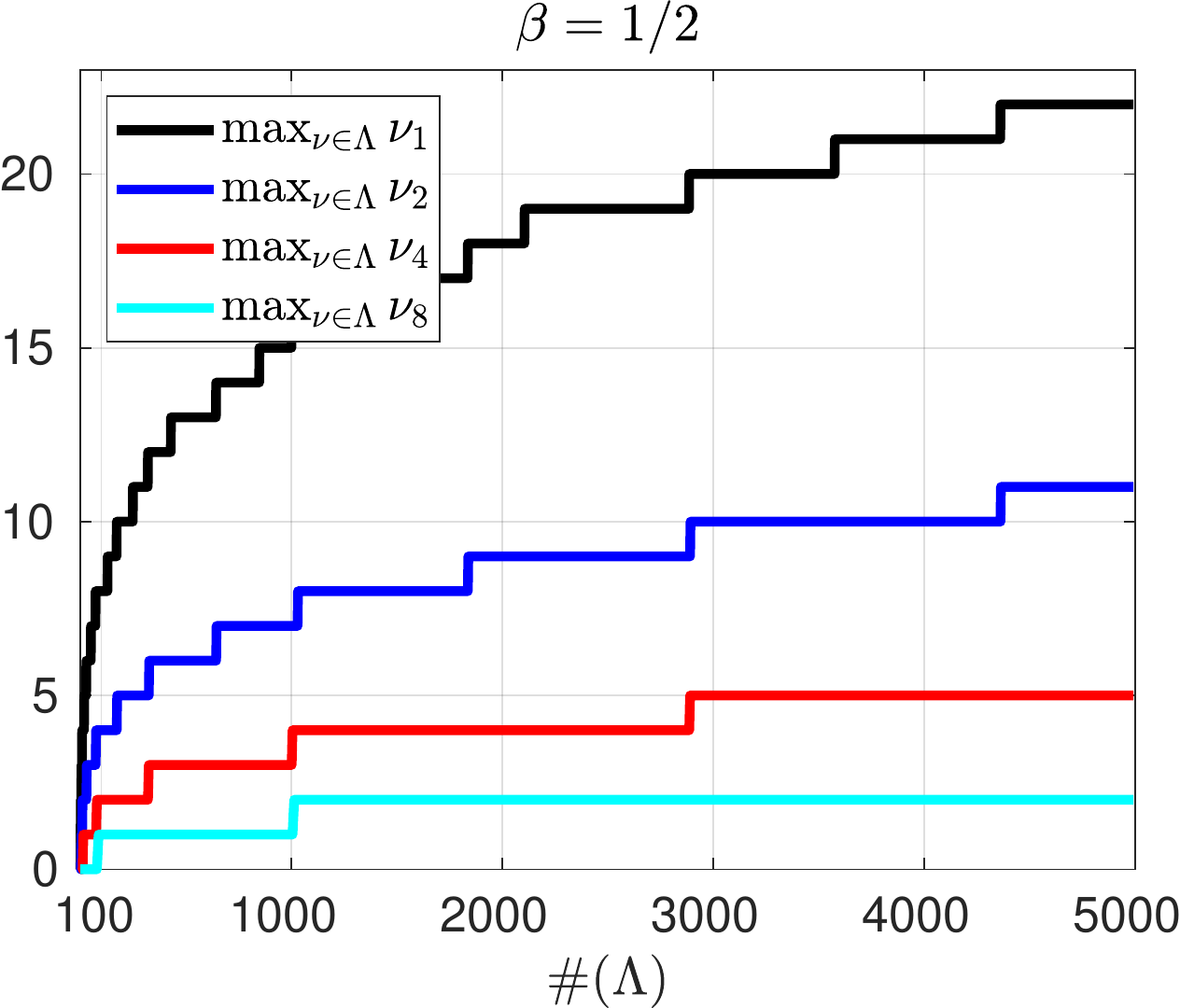} 
\caption{
Largest indices in $\Lambda$ associated to $\psi_{0,0}$, $\psi_{1,0}$, $\psi_{2,0}$, $\psi_{3,0}$:  
$\beta=\frac18$ (left), $\beta=\frac14$ (center), $\beta=\frac12$ (right).
}
\label{fig:activations}
\end{figure}

\begin{figure}[!h]
\center
\includegraphics[scale=0.4]{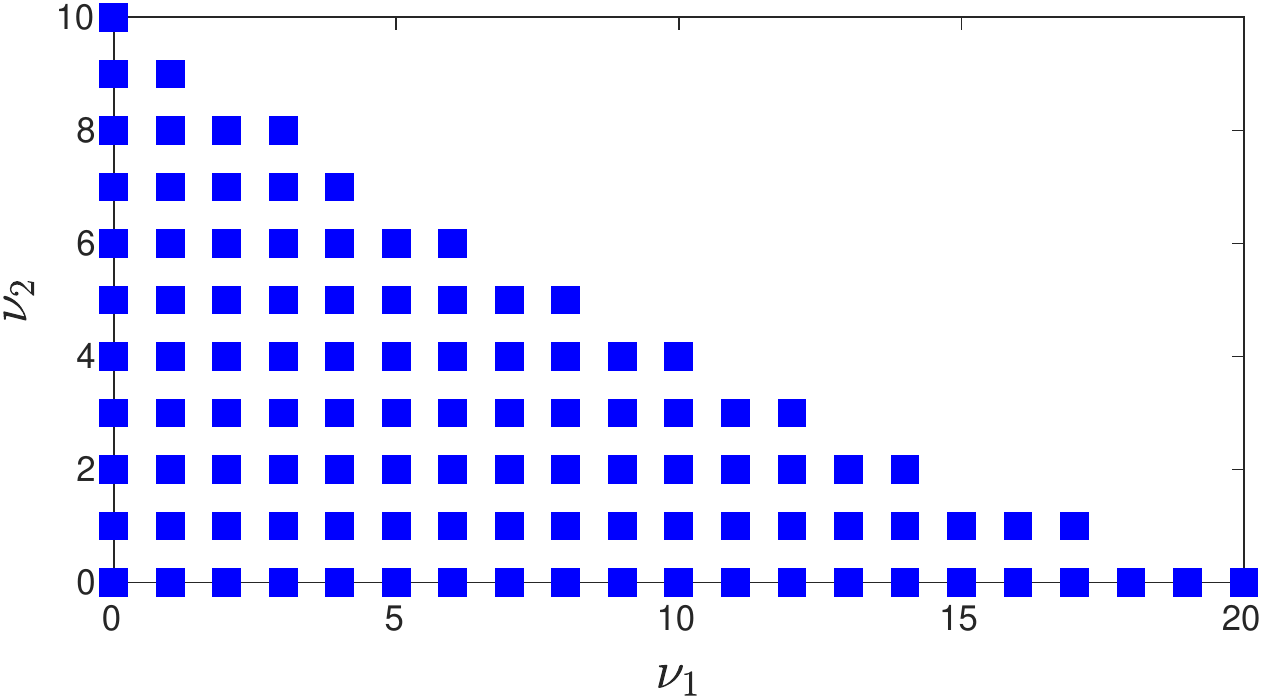} 
\includegraphics[scale=0.4]{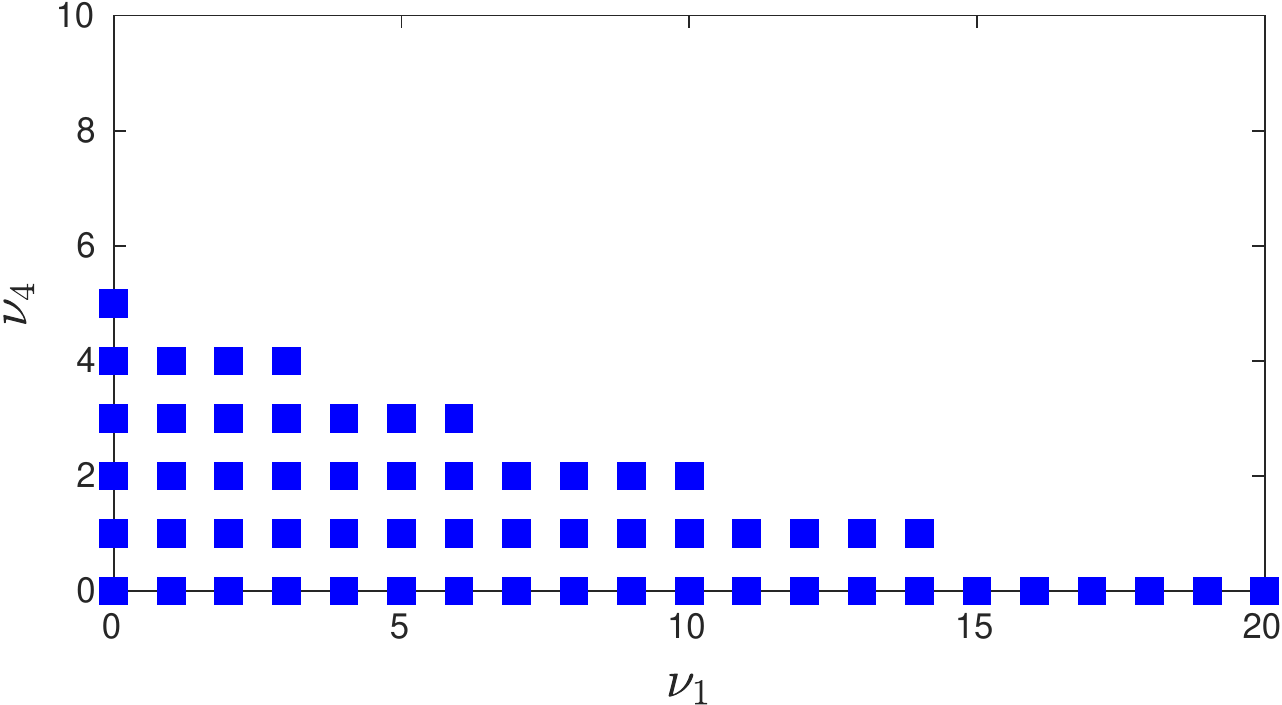} 
\includegraphics[scale=0.4]{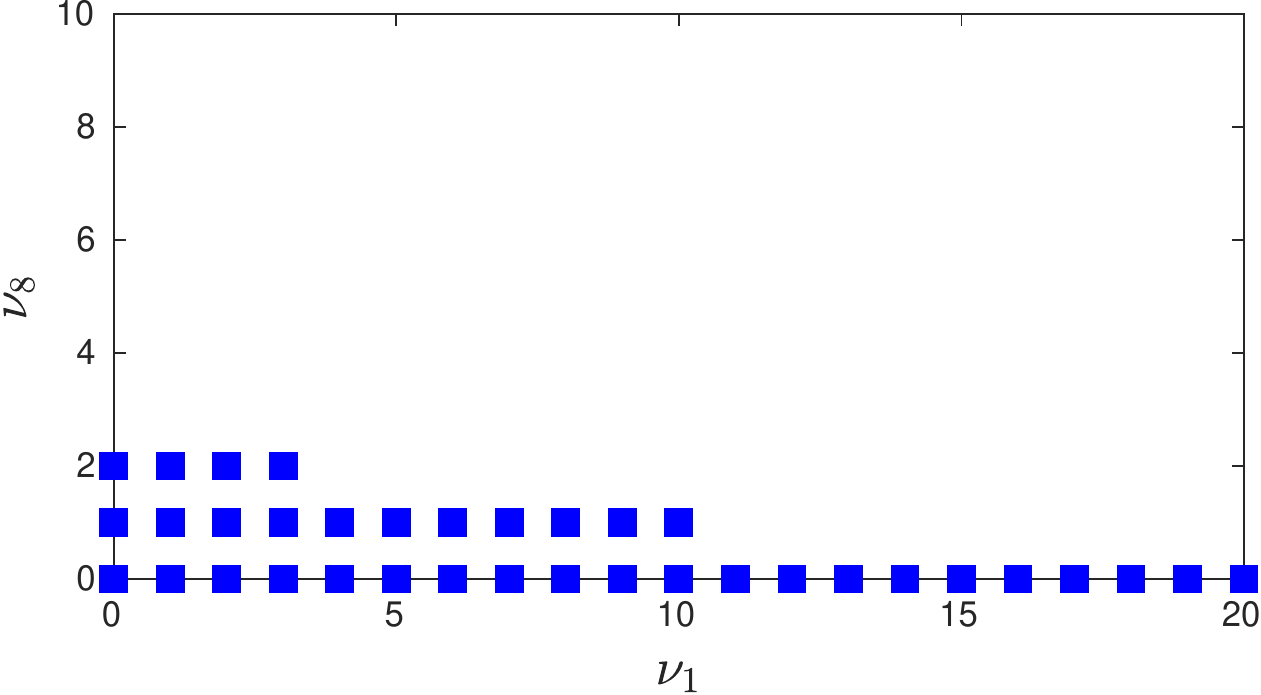}\\
\includegraphics[scale=0.4]{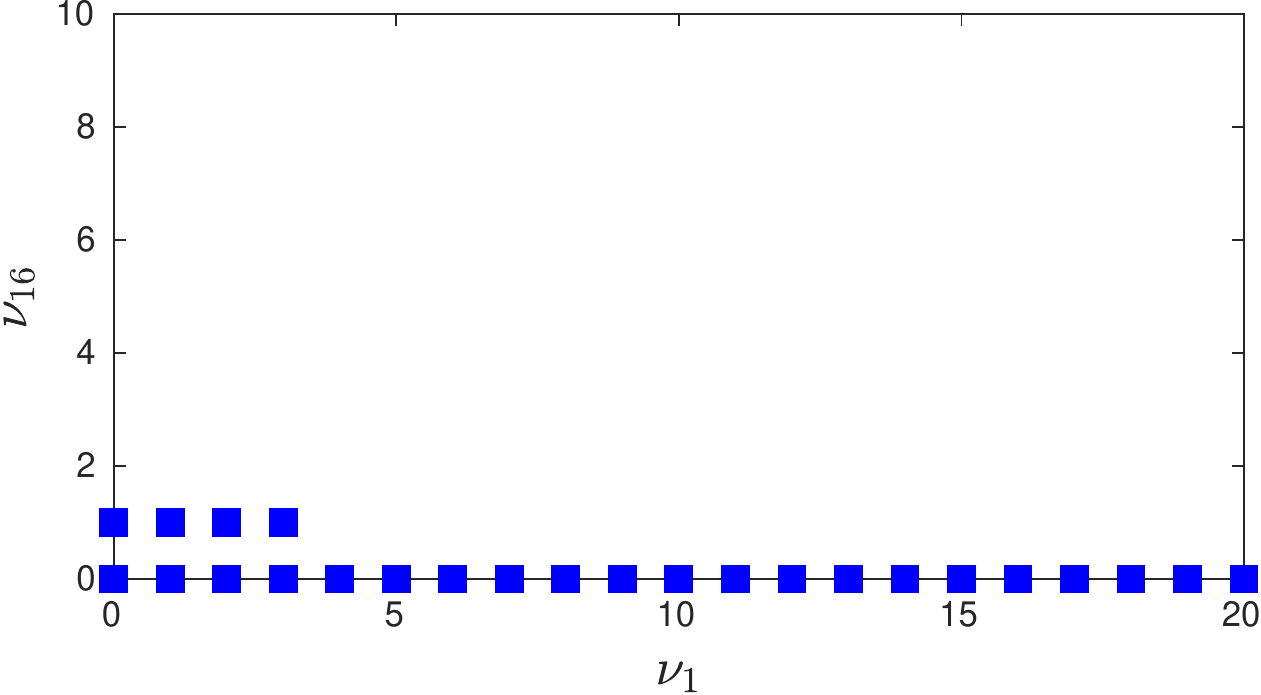} 
\includegraphics[scale=0.4]{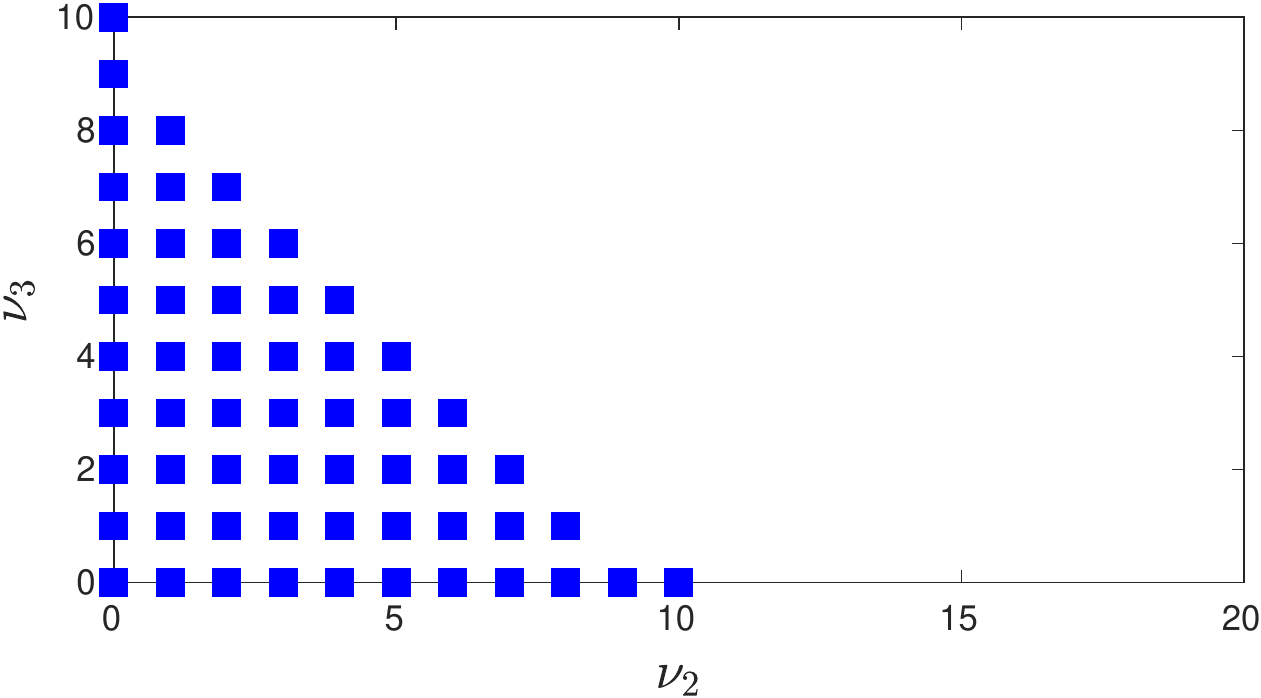} 
\includegraphics[scale=0.4]{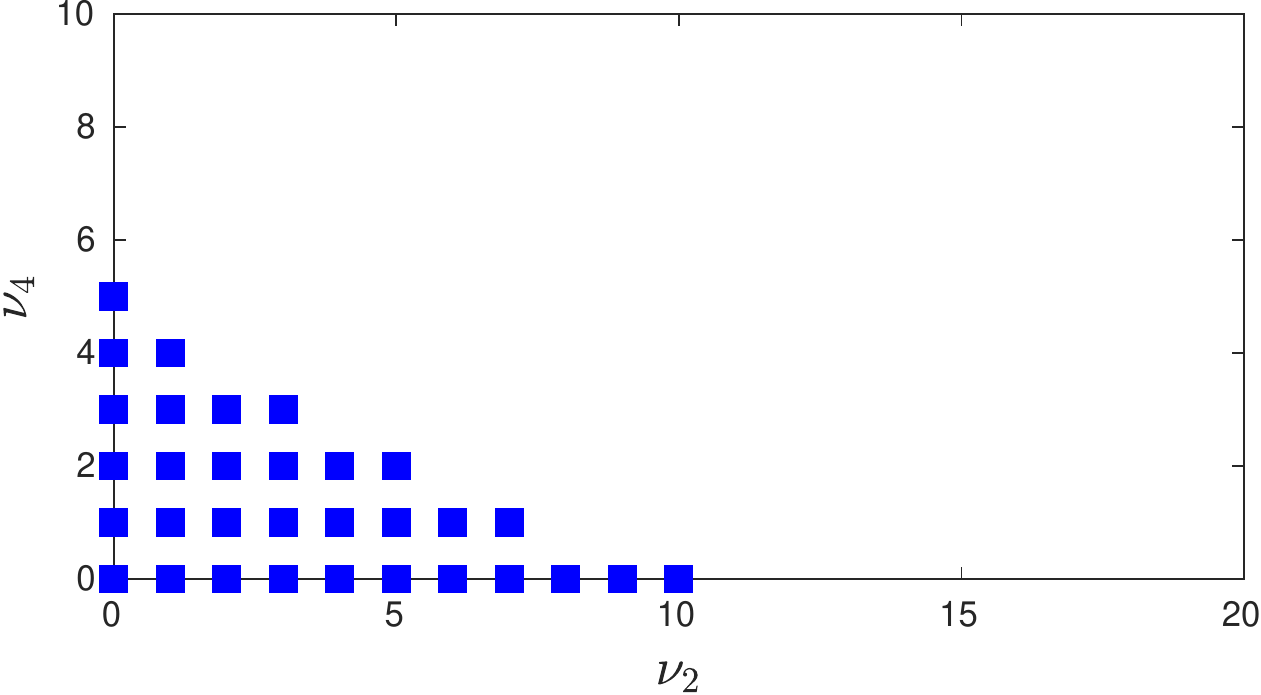} 
\caption{Some sections of the set $\Lambda_n$ with $n=4000$ and $\beta=\frac12$. 
}
\label{fig:sections_beta05}
\end{figure}

Figure \ref{fig:activations} shows the growth of  
the largest (component of the) indices of the Hermite polynomials associated to the 
random variables $y_{0,0}$, $y_{1,0}$, $y_{2,0}$, $y_{3,0}$ that multiply the 
first four levels of the Schauder basis 
$\psi_{0,0}$, $\psi_{1,0}$, $\psi_{2,0}$, $\psi_{3,0}$  
in \eqref{eq:schauder_basis}. 
An example of 
the interactions between different levels of the Schauder basis
is shown 
in Figure~\ref{fig:sections_beta05} when $\beta=\frac12$. 
Notice that anisotropy is between components associated to different levels of the Schauder basis, \emph{e.g.}~between $\psi_{0,0}$ and $\psi_{1,0}$ as shown in 
Figure~\ref{fig:sections_beta05}-top-left,
 but not between components associated to the same level, 
\emph{e.g.}~between $\psi_{1,0}$ and $\psi_{1,1}$ as in Figure~\ref{fig:sections_beta05}-bottom-center, whose sections remain isotropic regardless of the value of $\beta$. 

For the numerical validation of the estimate 
\eqref{eq:est_wls_bb} 
we approximate 
\begin{equation}
\label{eq:mc_exp_ref}
\mathbb{E}(\| u_h - u_{h,n}^C \|_{\cV_2} )
\approx 
\mathbb{E}(\| u_{h,{\textrm{ref}}}^C - u_{h,n}^C \|_{\cV_2} )
\end{equation}
by replacing $u_h$ with a 
reference solution $u_{h,{\textrm{ref}}}^C= \sum_{\nu \in \Lambda_{{\textrm{ref}}}} u_{h,\textrm{ref},\nu}^C H_\nu$, 
which is computed as a weighted least-squares estimator of $u_h$ with 
$n_{\textrm{ref}}:=\#(\Lambda_{\textrm{ref}})
=
5 \times 
10^3$  
and 
$m_{\textrm{ref}} = 20 \, n_{\textrm{ref}} \, \lceil \ln n_{\textrm{ref}} \rceil$
random samples. 
The expectation 
in the right-hand side of 
\eqref{eq:mc_exp_ref}
is 
estimated 
as a Monte Carlo average with $\widetilde{m}$ runs.

\begin{figure}[!h]
\center
\includegraphics[scale=0.8]{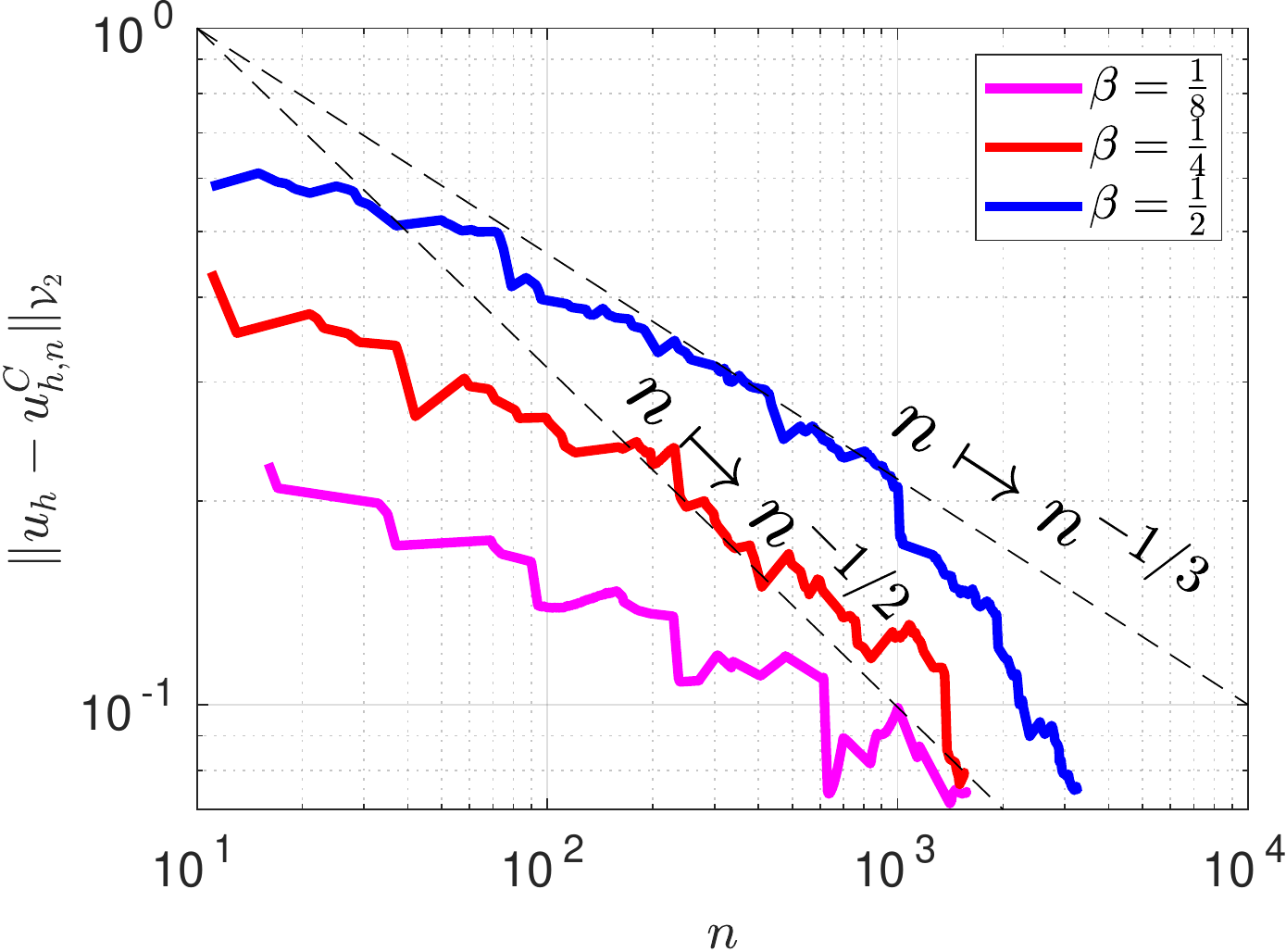} 
\caption{Monte Carlo averages of 
$\| u_{h,{\textrm{ref}}}^C - u_{h,n}^C \|_{\cV_2}$ 
with $\widetilde{m}=5$: 
$\beta=\frac18, \frac14, \frac12$, 
$\tau=1$.
}
\label{fig:errors_ref_1}
\end{figure}

Figure~\ref{fig:errors_ref_1} shows the convergence rates of the error \eqref{eq:mc_exp_ref} estimated with Monte Carlo and $\widetilde{m}=5$. 
The convergence plot depends on $\beta$.
Low values of $\beta$ like $\beta=\frac18$ 
produce a staircase convergence plot, 
with a substantial reduction of the error when a relevant group of important indices is activated. 
With higher values of $\beta$ like $\beta=\frac14$ or $\beta=\frac12$  
the convergence plot approaches the optimal asyptotic rate $n^{-1/2}$ predicted by \eqref{eq:est_wls_bb}. 

\begin{thebibliography}{77}

\bibliographystyle{plain}




\bibitem{AT03}
Ayache, A. and Taqqu, M. S.,
\newblock {\em Rate optimality of wavelet series approximations of fractional Brownian motion}.
\newblock  J. Fourier Anal. Appl. 9(5):451--471, 2003.

\bibitem{BCDM} 
Bachmayr, M., Cohen, A., DeVore, R. Migliorati, G.,
\newblock {\em Sparse polynomial approximation of parametric elliptic PDEs. Part II: lognormal coefficients}.
\newblock ESAIM:M2AN, 51:341--363, 2017.

\bibitem{BCM} 
Bachmayr, M., Cohen, A. and Migliorati, G., 
\newblock {\em Representations of Gaussian random fields and approximation of elliptic PDEs with lognormal coefficients}.
\newblock J. Fourier Anal. Appl. 1--29, 2017

\bibitem{BL}
Bergh, J. and L\"ofstr\"om, J.,
\newblock {\em Interpolation spaces, an introduction}.
\newblock  Springer, 1976



\bibitem{B98}  
Bogachev, V.I.,
\newblock {\em Gaussian measures}. 
\newblock American Mathematical Society Providence. Vol 62, 1998.

\bibitem{Cha} 
Charrier J., 
\newblock {\em Strong and weak error estimates for elliptic partial differential equations with random coefficients}.
\newblock  SIAM J. Numer. Anal. 50(1), 216--246, 2012.

\bibitem{Cia} 
Ciarlet, P.,
\newblock {\em The Finite Element Method for Elliptic Problems}.
\newblock North Holland Publ. 1978.

\bibitem{C}
Cohen, A.,
\newblock {\em Numerical analysis of wavelet methods}.
\newblock Elsevier, 2003.

\bibitem{CD} 
Cohen, A., DeVore, R., 
\newblock {\em Approximation of high dimensional parametric PDEs}.
\newblock Acta Numer., 24:1--159, 2015.


\bibitem{CM2016} 
Cohen, A., Migliorati, G.,
\newblock {\em Optimal weighted least-squares methods}. 
\newblock SMAI Journal of Computational Mathematics, 3:181--203, 2017. 

\bibitem{CM2016b} 
Cohen, A., Migliorati, G.,
\newblock {\em Multivariate Approximation in Downward Closed Polynomial Spaces}. 
\newblock 
Contemporary Computational Mathematics - A celebration of the 80th birthday of Ian Sloan, Springer, 2018.  


\bibitem{CrLea67}
Cramer, H., and Leadbetter, M. R., 
\newblock {\em Stationary and related stochastic processes: Sample function properties and their applications}. 
\newblock  Wiley, New York, 1967.

\bibitem{DS17}
Dashti M., Stuart A. M.,
\newblock {\em The Bayesian Approach to Inverse Problems}. 
\newblock Handbook of Uncertainty Quantification, Springer, 2015.

\bibitem{GS2009}
Galvis, J., Sarkis, M.,
\newblock {\em Approximating infinity-dimensional stochastic Darcy's equations without uniform ellipticity}.   
\newblock SIAM J. Numer. Anal. 47:3624--3651, 2009.


\bibitem{Gitt10}
Gittelson, C.,
\newblock {\em Stochastic Galerkin discretization of the log-normal isotropic diffusion problem}.  
\newblock Mathematical Models and Methods in Applied Sciences, 20(02):237--263, 2010.

\bibitem{G11}
Grisvard, P.,  
\newblock {\em Elliptic problems in nonsmooth domains}. 
\newblock Society for Industrial and Applied Mathematics, 2011.

\bibitem{Hac}
Hackbusch, W., 
\newblock {\em Elliptic differential equations - theory and numerical treatment}.  
\newblock Springer, 1992.

\bibitem{HS2014}
Hoang, V.H., Schwab, C.,
\newblock {\em N-term Galerkin Wiener chaos approximation rates for elliptic PDEs with lognormal Gaussian random inputs}.  
\newblock M3AS 24:797--826, 2014. 







\bibitem{RY13}
Revuz, D. and Yor, M., 
\newblock {\em Continuous martingales and Brownian motion}. 
\newblock Springer Science and Business Media, 2013.



\bibitem{Sz39}
Szeg\"o G.,
\newblock {\em Orthogonal polynomials}.
\newblock American Mathematical Society, Vol. 23, 1939.

\end{thebibliography}
\end{document}